\documentclass[a4paper,reqno]{amsart}
\usepackage{amsmath,amscd}
\usepackage{amssymb,stmaryrd}
\usepackage[french]{babel}
\usepackage[utf8]{inputenc}
\usepackage{hyperref}
\usepackage{enumerate}
\newcommand{\bbar}{\overline}
\newcommand{\cqfd}{\hfill$\Box$}

\theoremstyle{plain}
\newtheorem{thm}{Théorème}[section] 
\newtheorem{defn}[thm]{Définition} 
\newtheorem{lem}[thm]{Lemme}
\newtheorem{prop}[thm]{Proposition}
\newtheorem{princ}[thm]{Principe}
\newtheorem{coro}[thm]{Corollaire}

\theoremstyle{definition}
\newtheorem{rem}[thm]{Remarque}
\begin{document}

\title{Th\'eor\`eme d'Eilenberg-Zilber en
homologie cyclique enti\`ere}

\author{A. Bauval}\thanks{Anne Bauval, bauval@math.univ-toulouse.fr, Laboratoire \'Emile Picard (UMR5580), Universit\'e Toulouse~III}

\date{Janvier 1998}

\begin{abstract}Pour des modules simpliciaux, le th\'eor\`eme d'Eilenberg-Zilber
classique \'enonce l'existence d'un produit $sh : M\otimes N\to
M\times N$ (le shuffle) et d'un coproduit $AW : M\times N\to
M\otimes N$ (l'application d'Alexander-Whitney), quasi-inverses. Une version cyclique de ce th\'eor\`eme a été établie (\cite{HJ}), prouvant l'existence de ``coextensions'' $sh_\infty$ et $AW_\infty$ de $sh$ et $AW$, par une m\'ethode de
mod\`eles acycliques. Par ailleurs, une formule explicite pour
$sh_\infty$ a \'et\'e d\'ecouverte par divers auteurs.

Nous résolvons le problème restant : expliciter de m\^eme $AW_\infty$, ainsi que les
homotopies par lesquelles $sh_\infty$ et $AW_\infty$ sont
quasi-inverses et quasi-(co)-associatifs, puis montrons
que toutes les applications explicit\'ees s'\'etendent
contin\^ument aux complexes cycliques entiers (associ\'es \`a des alg\`ebres
norm\'ees).
\end{abstract}

\maketitle

\noindent{\bf Introduction.--} Pour des modules simpliciaux, le th\'eor\`eme
d'Eilenberg-Zilber classique \'enonce l'existence d'un produit $sh : M\otimes
N\to M\times N$ (le shuffle) et d'un coproduit $AW : M\times N\to
M\otimes N$ (l'application d'Alexander-Whitney), quasi-inverses.

Une version cyclique de ce th\'eor\`eme a \'et\'e \'etablie par
Hood-Jones \cite{HJ}~: ils prouvent l'existence de
\og~coextensions\fg\  $sh_\infty$ et $AW_\infty$ de $sh$ et $AW$,
par une m\'ethode de mod\`eles acycliques. Par ailleurs, une formule
explicite pour co\'etendre $\bbar{sh}$ (la version normalis\'ee de
$sh$)  a \'et\'e d\'ecouverte par divers auteurs (\cite{L},
bibliographie du chapitre 4)~:
$\bbar{sh}_\infty$ est simplement constitu\'e du shuffle
$\bbar{sh}$ et du shuffle cyclique $\bbar{sh}'$.

Mais la question d'expliciter de m\^eme $\bbar{AW}_\infty$, ainsi
que les homotopies par lesquelles $\bbar{sh}_\infty$ et
$\bbar{AW}_\infty$ sont quasi-inverses et quasi-(co)-associatifs
restait ouverte. Or elle se pose de fa\c con cruciale lorsqu'on
veut passer \`a l'homologie cyclique enti\`ere ou asymptotique,
donc v\'erifier la continuit\'e de ces applications. Nous
fournissons une solution compl\`ete \`a ce probl\`eme
(\og~red\'ecouvrant\fg\  au passage le shuffle cyclique), puis
montrons que toutes les applications explicit\'ees s'\'etendent
contin\^ument aux complexes cycliques entiers (associ\'es \`a des
alg\`ebres norm\'ees $A$ et $B$), pour donner un quasi-isomorphisme
$$\Omega_\varepsilon (A)\otimes\Omega_\varepsilon
(B)\simeq\Omega_\varepsilon (A\otimes B).$$

Cette solution est une \'etape vers une reformulation simplifi\'ee
de l'homologie cyclique asymptotique de Puschnigg (\cite{P1}), dont
la construction d'un produit et d'un coproduit (\cite{P2}) est
radicalement diff\'erente.

\noindent{\bf Plan.--} Les six sections de l'article sont, dans la mesure du
possible, ind\'ependantes, de mani\`ere \`a permettre au lecteur press\'e de
se limiter, en premi\`ere lecture, aux trois sections essentielles qui
sont~:
\begin{itemize}
\item{\ref{section:rappels}}.-- Rappels (principalement I.1 et I.2), o\`u l'on reformule
le probl\`eme et le lemme de perturbation de Brown, et les solutions
(partielles) ant\'erieures.
\item{\ref{section:coext}}.-- Coextensions
explicites de $\bbar{sh}$ et $\bbar{AW}$, o\`u l'on utilise ce lemme de
perturbation pour construire $\bbar{sh}_\infty$, $\bbar{AW}_\infty$, et
l'homotopie prouvant qu'ils sont quasi-inverses (on montre ensuite qu'ils sont
quasi-(co)-associatifs mais non quasi-(co)-commutatifs).
\item{\ref{section:entiere}}.-- Produit et coproduit en homologie cyclique enti\`ere,
o\`u l'on prouve la continuit\'e des applications ci-dessus.
\end{itemize}

Dans les sections \og~facultatives\fg\ on r\'esout (sections \ref{section:coext-long-2} et
\ref{section:denormalisation}) le m\^eme probl\`eme de coextensions (si possible de longueur 2)
dans le cadre non normalis\'e, et on montre (section \ref{section:genericite}) comment
\'etendre \`a tous les modules cycliques les
propri\'et\'es obtenues dans le cadre des modules cycliques
d'alg\`ebres.

\noindent{\bf Notations.--} $k$ est un anneau commutatif unitaire. Les ($k$-)
alg\`ebres seront toujours suppos\'ees unitaires. Pour tout
ensemble $E$, on notera $k[E]$ le $k$-module libre sur $E$. Pour
tout entier naturel $p$, $\Lambda [p]$ d\'esignera l'ensemble
cyclique $\Lambda ^{op}(p,-)$, et $\lambda [p]=k[\Lambda [p]]$ le
$k$-module cyclique engendr\'e. La $k$-alg\`ebre libre sur un
ensemble \`a $n+1$ \'el\'ements $\{a_0,\ldots, a_n\}$ sera not\'ee
$T_n$ (on peut aussi la d\'efinir comme l'alg\`ebre tensorielle
$T(k^{n+1})$). Le normalis\'e d'un module simplicial
(\'eventuellement cyclique) $P_*$ sera not\'e $\bbar P_*$ (\cite{L}
1.6.4). En particulier pour $P_*=C_*(A)$ (associ\'e \`a une
alg\`ebre $A$), $\bar C_n(A)=A\otimes (\bar A)^{\otimes
n}=\Omega^n(A)$. Dans $\Omega(A\otimes B)$, on notera abusivement
$a$ au lieu de $a\otimes 1$ et $b$ au lieu de $1\otimes b$ (pour
$a\in A$ et $b\in B$). Les op\'erateurs usuels $\bar b$ et $\bar B$
(sur les complexes normalis\'es) seront not\'es simplement $b$ et
$B$.

\section{Rappels}\label{section:rappels}

\subsection{Formulation du probl\`eme}\label{subsection:formulation}
Soient $P$ et $Q$ deux modules simpliciaux. On note
$P_*,Q_*$ les complexes associ\'es (munis de la diff\'erentielle de
Hochschild $b:P_n\to P_{n-1}$). Les applications
$$\begin{matrix}AW & : & (P\times Q)_* &\to &P_*\otimes Q_*\\
sh:&:&P_*\otimes Q_*&\to &(P\times Q)_* \end{matrix}$$
 sont des quasi-isomorphismes naturels, inverses l'un de
l'autre, et $sh$ est associative et (avec la graduation)
commutative. De plus, ces deux applications passent aux
normalis\'es, et d\'efinissent des quasi-isomorphismes $\bbar {AW},
\bbar{sh}$ entre $\bbar{(P\times Q)_* }$ et\break $\bbar
{P_*}\otimes \bbar{Q_*}$ (\cite{L} 1.6.4--1.6.12). On peut aussi
remarquer que $AW$ est coassociative.

Pr\'ecisons ce qu'on entend par \og~$sh$ est (avec la graduation)
commutative\fg\ . On introduit deux op\'erateurs $\sigma $ (qui
commutent \`a $b$) en posant
$$\begin{matrix}\sigma &:&P\otimes Q&\to&Q\otimes P\\
&&x_p\otimes y_q&\mapsto&(-1)^{pq}y_q\otimes x_p\end{matrix}
\qquad\qquad\begin{matrix}\sigma &:&P\times Q&\to &Q\times P\\
&& x_n\otimes y_n&\mapsto&y_n\otimes x_n\end{matrix}$$
(pour $x_p\in
P_p$, $y_q\in Q_q$, $x_n\in
P_n$ et $y_n\in Q_n$), et l'on v\'erifie que $\sigma
sh=sh \sigma $.

Si $P,Q$ sont des modules non seulement simpliciaux mais
cycliques, $P_*, Q_*$ sont \'egalement munis de l'op\'erateur de
Connes $B:P_n\to P_{n+1}$ (qui passe aux normalis\'es). On pose
$C^-_d(P_*)=\prod_i P_{d+2i}$ et on munit $C^-_*(P_*)$ de la
diff\'erentielle $b+B:C^-_d(P_*)\to C^-_{d-1}(P_*)$. H\'elas, $AW$
et $sh$ ne commutent pas \`a $B$ (m\^eme sur les normalis\'es) et
m\^eme, les $(b,B)$-bicomplexes $(P\times Q)_*$ et $P_*\otimes
Q_*$ ne sont g\'en\'eralement pas quasi-isomorphes. Mais on peut,
pour prouver que $C^-((P\times Q)_*)$ et $C^-(P_*\otimes Q_*)$ le
sont n\'eanmoins, construire ce que Hood et Jones (\cite{HJ}) appellent
des {\sl coextensions\/} de $AW$ et $sh$ (et de $\bbar{AW}$ et
$\bbar{sh}$). On suppose dans
la d\'efinition suivante que $P_*, Q_*$ ($*\geq 0$) sont munis de
diff\'erentielles $b$ de degr\'e $-1$, $B$ de degr\'e $+1$, telles
que $bB=-Bb$.

\begin{defn}
Soit $f:P_*\to Q_*$ un
morphisme de $b$-complexes. Une coextension de $f$ est une suite
$f_\infty=(f_k)_{k\geq 0}$, avec $f_0=f$, $f_k:P_*\to Q_{*+2k}$,
telle que $[b+B,f_\infty]=0$, c'est-\`a-dire $$ [b,
f_k]+[B,f_{k-1}]=0$$ (par convention, $f_{-1}=0$). Cette
coextension sera dite de longueur $2$ si $f_k=0$ pour $k\geq 2$.
\end{defn}

\subsection{Lemme de perturbation et normalisation}\label{subsection:perturbnorma}

Avant de rappeler les m\'ethodes de Hood-Jones et de Kassel pour
construire de telles coextensions, mentionnons (en le
particularisant \`a notre contexte) un lemme de Brown, remis \`a
l'honneur par Kassel dans le cadre cyclique (\cite{K1} lemme 1.3 et
proposition 7.1), et dont une premi\`ere application imm\'ediate
donne un quasi-isomorphisme naturel entre $C^-(Q_*)$ et
$C^-(\bbar{Q_*})$ \cite{K1}, application 7.3.a). On adopte les m\^emes
hypoth\`eses et notations que
dans la d\'efinition I.1.

\begin{lem}\label{lem:perturbation} Soient $f:P_*\to Q_*, g:Q_*\to P_*$, morphismes de
$b$-complexes tels que $gf=1$ et $\varphi :Q_*\to Q_{*+1}$ une
homotopie telle que $fg=1+b\varphi +\varphi b$. On pose
$$\begin{matrix} \varphi _\infty&=&\varphi \sum_{m\geq 0}(B\varphi )^m
&=&\sum_{m\geq 0}(\varphi B)^m\varphi,\\
f_\infty&=&(1+\varphi _\infty B)f,&&\\
g_\infty&=&g(1+B\varphi _\infty),&&\\
b_\infty&=&b+gB(1+\varphi _\infty B)f&=&b+g(1+B\varphi
_\infty)Bf.\end{matrix}$$
Alors,
$$\begin{matrix}b_\infty^2&=&0,\\
f_\infty b_\infty&=&(b+B)f_\infty,\\
g_\infty (b+B)&=&b_\infty g_\infty,\\
f_\infty g_\infty&=&1+(b+B)\varphi _\infty+\varphi
_\infty(b+B).\end{matrix}$$
Si de plus $\varphi $ est \og~sp\'eciale\fg~, c'est-\`a-dire
v\'erifie $$ \varphi f=0,\qquad g\varphi =0,\qquad \varphi ^2=0,$$
alors $g_\infty f_\infty=1$ et $\varphi_\infty$ est
sp\'eciale.
\end{lem}

($\varphi $ peut toujours \^etre rendue \og~sp\'eciale~\fg, \cite{K1} 
remarque 1.2). Le quasi-iso\-morphisme (naturel) entre $C^-(Q_*)$ et
$C^-(\bbar{Q_*})$ s'obtient en appliquant ce lemme \`a
$P_*=\bbar{Q_*}$ et $g=j:Q_*\to \bbar{Q_*}$ la projection
canonique. Il existe en effet $i:\bbar{Q_*}\to Q_*$ et $\varphi $
v\'erifiant les hypoth\`eses du lemme. Or $j$ commute \`a $B$ (par
d\'efinition de $B$ sur $\bbar{Q_*}$), d'o\`u $b_\infty=b+B$ (et
$j_\infty=j$).

D'apr\`es ce r\'esultat, si des coextensions $AW_\infty,
sh_\infty$ existent, on peut en d\'e\-duire des coextensions
$\bbar{AW}_\infty:=(j_\infty\otimes j_\infty)AW_\infty i_\infty$ et
$\bbar{sh}_\infty:=j_\infty sh_\infty(i_\infty\otimes i_\infty) $
(on aura bien  $\bbar{AW}_0=(j\otimes j)AW
i=\bbar{AW}ji=\bbar{AW}$, et $\bbar{sh}_0=j sh (i\otimes
i)=\bbar{sh}(j\otimes j)(i\otimes i)=\bbar{sh}$. On pourrait croire
(\cite{HJ}, proof of theorem 2.3.a) que la r\'eciproque est aussi
simple~: si des coextensions $\bbar{AW}_\infty, \bbar{sh}_\infty$
existent, on a, certes (\cite{GJ1}), un quasi-isomorphisme
$g_\infty:=(i_\infty\otimes i_\infty)\bbar{AW}_\infty j_\infty$,
$f_\infty:=i_\infty\bbar{sh}_\infty(j_\infty\otimes j_\infty)$,
mais ces $g_\infty$ et $f_\infty$ ne sont que des
coextensions de $g_0=(ij\otimes ij)AW\neq AW$ et de $f_0=ijsh
\neq sh$ (on a donc seulement $(j\otimes j)g_0=(j\otimes j)AW$ et
$jf_0=jsh$). Ce probl\`eme sera r\'esolu au \S \ref{section:denormalisation}.

D'autre part, le lemme ci-dessus sera r\'eutilis\'e au \S \ref{section:coext}.

\subsection{M\'ethode de Hood-Jones}\label{subsection:Hood-Jones}

J'exposerai uniquement la m\'ethode de \cite{HJ} pour co-\'etendre
un $f_0:\bbar{P_*}\otimes\bbar{Q_*}\to\bbar{(P\times Q)_*}$, car la
solution pour un $g_0: \bbar{(P\times
Q)_*}\to\bbar{P_*}\otimes\bbar{Q_*}$ est analogue, et car il faut
un argument suppl\'ementaire (cf. \S\ ci-dessus et \S \ref{section:denormalisation})
pour passer aux non-normalis\'es. Soit donc
$f_0:\bbar{P_*}\otimes\bbar{Q_*}\to\bbar{(P\times Q)_*}$, un
quasi-isomorhisme (naturel) de $b$-complexes, \'egal \`a ${\rm id}$
en degr\'e $0$. (Ou plus g\'en\'eralement~: un morphisme homotope
\`a $\bbar{sh}$).

On construit les $f_k$ par
r\'ecurrence sur $k$ puis par sous-r\'ecurrence sur le degr\'e
total $d=p+q$ de l'\'el\'ement $\bar x\otimes \bar y\in\bbar
P_p\otimes \bbar Q_q$ auquel $f_k$ s'applique. Puisqu'on veut $f_k$
naturel, on va d\'efinir $f_k(\bar x\otimes \bar y)$ uniquement pour
$x=1_p\in\lambda [p], y=1_q\in\lambda [q]$,
puis \'etendre $f_k$ \`a $P,Q,x,y$ quelconques en utilisant
l'universalit\'e des $\lambda [n]$ (\cite{HJ} lemme 2.1), qui n'est en
fait qu'un avatar du lemme de Yoneda.

L\`a encore, un probl\`eme se pose (\`a moi) dans la m\'ethode de
\cite{HJ}~: si $\bbar z=f_k(\bbar{1_p}\otimes \bbar{1_q})$, et si
$i_x:\lambda [p]\to P, i_y:\lambda [q]\to Q$ sont les morphismes
canoniques de modules cycliques tels que $i_x(1_p)=x$ et
$i_y(1_q)=y$, on voudrait que
$$f_k(\bar x\otimes\bar y)=\bbar{(i_x\times i_y)(z)}.$$
Or rien a priori n'assure que cette d\'efinition soit possible~:
il faudrait pour cela que lorsque $x$ ou $y$ est d\'eg\'en\'er\'e,
$(i_x\times i_y)(z)$ le soit aussi. Nous verrons (remarque \ref{rem:normalisable})
qu'en fait cette condition est automatiquement r\'ealis\'ee pour
$k=1$, et que pour $k\geq 2$ on a m\^eme $f_k=0$. Par contre, le
probl\`eme se pose r\'eellement dans l'autre sens, lorsqu'on veut
co\'etendre un $g_0:
\bbar{(P\times Q)_*}\to\bbar{P_*}\otimes\bbar{Q_*}$.

Pour \'eviter
les deux probl\`emes ci-dessus, reprenons la m\'ethode
de \cite{HJ} mais sans passer par les normalis\'es, et construisons
directement une coextension $f_\infty$ d'un
$f_0: P_*\otimes Q_*\to (P\times Q)_*$. L'\'etape de
r\'ecurrence, c'est-\`a-dire la construction d'un $z\in k[(\Lambda
[p]\times\Lambda [q])_{d+2k}]$ tel que
$b(z)=Z:=(f_kb-[B,f_{k-1}])(1_p\otimes 1_q)$ se fait en remarquant
que (par hypoth\`ese de r\'ecurrence) $bZ=0$ donc (par acyclicit\'e
de $(\lambda [n],b)$ en degr\'es $\geq 2$) $z$ existe d\`es que
$d+2k-1\geq 3$. Il suffit donc d'initialiser la r\'ecurrence, en
construisant \og~\`a la main\fg\  $f_k(1_p\otimes 1_q)$ pour $k=1$
et $d\leq 1$. Dans leur contexte (normalis\'e), Hood-Jones exhibent
une telle initialisation pour $\bbar {sh}_1$, puis en d\'eduisent
une initialisation pour un $f_0$ \og~quelconque\fg~, c'est-\`a-dire
de la forme $\bbar {sh}+bh+hb$, en posant $f_1=\bbar{sh}_1+Bh+hB$.
Pour adapter cela au cadre non normalis\'e, il suffit d'initialiser
$sh_1$, puis d'appliquer leur argument pour $f_0=sh+bh+hb$.

Nous donnerons au \S \ref{section:denormalisation} (par une autre m\'ethode) une preuve de
l'existence de $sh_\infty$, ce qui prouve que l'initialisation de
$sh_1$ est possible (remarque \ref{rem:V.5}).

Remarquons qu'une fois $sh_\infty$ construit
(par la m\'ethode ci-dessus, ou par celle du \S V), il est en
r\'ealit\'e inutile, pour co-\'etendre $f_0$, d'initialiser $f_1$
puis d'appliquer (de nouveau) la m\'ethode des mod\`eles
acycliques~: il suffit de poser $f_1=sh_1+Bh+hB$ et $f_k=sh_k$ pour
$k\geq 2$. Dans le contexte normalis\'e, le raisonnement serait le
m\^eme pour co\'etendre $\bbar {sh}+bh+hb$, \`a partir d'une
coextension de $\bbar{sh}$ que nous fournirons au \S IV.

\subsection{M\'ethode de Kassel}\label{subsection:Kassel}

Quelques mois auparavant, dans le cadre a priori plus restreint
(mais voir \S \ref{section:genericite}) des modules cycliques associ\'es \`a des
alg\`ebres, et pour co\'etendre seulement $\bbar{sh}$ (et
pas $\bbar{AW}$), Kassel \cite{K2} utilisait une m\'ethode voisine. Les
\og~mod\`eles acycliques\fg\  $\lambda [n]$ sont remplac\'es, dans
ce contexte, par les modules cycliques $C_*(T_n)$.

Kassel cherchait $\bbar{sh}_k$ sous la forme
$F^{(k)}\circ\bbar{sh}$ (ce qui n'est pas restrictif puisque
$\bbar{AW}\ \bbar{sh}=1$), les $F^{(k)}$ \'etant constitu\'es des
$F^{(k)}_d:\bar C_d(A\otimes B)\to\bar C_{d+2k}(A\otimes B)$.
L'\'etape de r\'ecurrence (sur $k$ et $d$) consiste alors \`a
trouver $$F^{(k)}_d((a_0\otimes b_0){\rm d}(a_1\otimes
b_1)\ldots {\rm d}(a_d\otimes b_d))=z\in \bar C_{d+2k}(T_d\otimes
T_d)$$ tel que $bz=Z$, o\`u $Z$ est (par hypoth\`ese de
r\'ecurrence) un certain cycle dans $\bar C_{d+2k-1}(T_d\otimes
T_d)$. Par acyclicit\'e des $C_*(T_n)$ en degr\'es $\geq 2$, $z$
existe d\`es que $d+2k-1\geq 3$, et il suffit donc de construire
\`a la main l'initialisation $F^{(1)}_0$ et $F^{(1)}_1$. ($bZ$ est
bien nul seulement dans le complexe associ\'e \`a $T_d\otimes T_d$,
et non dans celui associ\'e \`a $T_{2d+1}$, dont l'alg\`ebre
pr\'ec\'edente n'est qu'un quotient~; \`a cause probablement de
cette confusion, Kassel n'explicite que $F^{(1)}_0$~; nous
compl\`eterons son initialisation dans la remarque \ref{rem:calculinit}).

Mais \`a nouveau (comme au \S \ref{subsection:Hood-Jones}), le probl\`eme de
compatibilit\'e avec la norma\-lisation se pose, pour
\'etendre $F^{(k)}_d$ \`a un \'el\'ement quelconque de degr\'e $d$~:
lorsqu'on consid\`ere $z$ comme une fonction multilin\'eaire des
variables $a_0, b_0,\ldots, a_d, b_d$ et qu'on remplace $a_i\otimes
b_i$ par $1\otimes 1$ pour un certain $i\geq 1$, rien ne garantit
que le r\'esutat soit nul (sauf pour $d<2k$~: cf. remarque \ref{rem:normalisable}).

Retenons cependant de cette m\'ethode l'id\'ee de chercher
$\bbar{sh}_\infty$ sous la forme $F\circ\bbar{sh}$ (cf. \S \ref{section:coext}), et
celle de restreindre (en apparence) le contexte aux modules
cycliques d'alg\`ebres (cf. \S \ref{section:genericite}).

\subsection{Shuffles cycliques}\label{subsection:shuffle}

Divers auteurs ont, ind\'ependamment, trouv\'e une coextension de
$\bbar{sh}$, explicite, et m\^eme de longueur 2 (cf. \cite{L},
bibliographical comments on chapter 4). D\'ecrivons donc leur
$\bbar{sh}_1$.

Un $(p,q)$-shuffle cyclique est une permutation $\sigma$ sur $p+q$
\'el\'ements, telle que $\sigma (a_1,\ldots,a_p,b_1,\ldots,b_q)$
s'obtienne en effectuant (ind\'ependamment) une permutation
circulaire sur les $a_i$ et une sur les $b_j$, puis (tout en
conservant l'ordre obtenu sur les $a_i$ et celui sur les $b_j$), en
m\'elangeant les deux suites de telle fa\c con que $a_1$ reste
avant $b_1$. Notons $S'_{p,q}$ l'ensemble des $(p,q)$-shuffles
cycliques.

Si l'on note,
comme dans \cite{L} 4.2.1, l'action \`a gauche du groupe sym\'etrique
$S_n$ sur $\bar C_n(A)$ par $\tau  \bullet a_0{\rm d}a_1\ldots {\rm
d}a_n=a_0{\rm d}a_{\tau  ^{-1}(1)}\ldots {\rm d}a_{\tau ^{-1}(n)}$,
et qu'on pose $$(a_0{\rm d}a_1\ldots {\rm d}a_p)\perp (b_0{\rm
d}b_1\ldots{\rm d}b_q)=\sum_{\sigma \in S'_{p,q}}\varepsilon
(\sigma )\quad \sigma ^{-1}\bullet (a_0\otimes b_0){\rm
d}a_1\ldots{\rm d}a_p{\rm d}b_1\ldots{\rm d}b_q,$$
alors $\bbar {sh}_1:\bar C_p(A)\otimes\bar C_q(B)\to\bar
C_{p+q+2}(A\otimes B)$ peut \^etre choisi \'egal \`a~:
$$\bbar{sh}'(x_p \otimes y_q):=(-1)^p({\rm d}x_p)\perp({\rm
d}y_q).$$
(Les deux formules ci-dessus rectifient la confusion entre $\sigma
$ et $\sigma ^{-1}$ et l'oubli du $(-1)^p$, dans les deux formules
correspondantes de \cite{L} 4.3.2).

Plus directement (et sans utiliser la notation $\bullet$)
$$\bbar{sh}'(x_p\otimes y_q)=(-1)^p\sum_{\sigma \in S'_{p+1,q+1}}
\varepsilon (\sigma )\qquad \sigma ({\rm d}x_p{\rm d}y_q).$$
Ce $\bbar{sh}'$ r\'eappara\^\i tra naturellement au \S \ref{section:coext}. D'autres
choix sont bien s\^ur possibles pour $\bbar{sh}_1$, mais les
\og~initialisations\fg\  de Hood-Jones et de Kassel correspondaient
\`a celui-ci (cf. remarque \ref{rem:calculinit}). Quels que soient ces choix,
$\bbar{sh}_k$ peut \^etre pris nul pour $k\geq 2$. Ce (petit)
miracle sera \'eclairci au \S III (remarque \ref{rem:normalisable}).

\subsection{R\'esultats de Puschnigg}

Puschnigg a, le premier, construit un produit et un
coproduit suffisamment explicites, en homologie cyclique
 p\'eriodique, pour en contr\^oler la continuit\'e et en
d\'eduire les m\^emes op\'erations en homologie cyclique
enti\`ere. (Tant qu'ils restaient dans le cadre purement
alg\'ebrique, les auteurs pr\'ec\'edents se souciaient peu du
caract\`ere explicite de $\bbar{sh}_\infty, \bbar{AW}_\infty$, et
m\^eme (\cite{K2}, \cite{L}) -- gr\^ace au lemme des 5 -- de l'existence de
$\bbar{AW}_\infty$).

Alors que les rappels des
paragraphes pr\'ec\'edents seront des ingr\'edients techniques du
pr\'esent travail, la contribution de Puschnigg en a \'et\'e
l'ingr\'edient essentiel~: motivation, espoir, et strat\'egie (a
contrario). En effet, dans son introduction \cite{P2},
il
\begin{itemize}
\item pr\'etend qu'il est impossible d'expliciter les
r\'esultats ant\'erieurs~: ``~It is known that a homotopy inverse to
the shuffle map exists but cannot be described explicitely~''~, et
souligne l'aspect crucial d'une telle explicitation pour passer \`a
l'homologie cyclique enti\`ere~;
\item se limite, partant de ce constat d'\'echec, \`a
travailler sur les complexes p\'erio\-diques~: ``~To overcome these
difficulties we want to develop product operations from a
completely different point of view, due to Cuntz and Quillen
$[\ldots]$ in a purely ${\Bbb Z}/2{\Bbb Z}$-graded context~''~;
\item explique
pourquoi sa strat\'egie donne fatalement un produit peu
calculable, bien qu'~\og~explicite\fg\ (com\-po\-s\'e
de six applications explicites)~: ``~the tensor product of two
[quasi-free] algebras will be of cohomological dimension two
$[\ldots]$ their periodic cyclic (co)-homology can still be
calculated by a small quotient of the periodic cyclic complex
$[\ldots]$. All calculations become much more elaborate however due
to the lack of a lifting property $[\ldots]$ the construction of a
product map will be considerably harder than that of the coproduct
map~''~.
\end{itemize}

Il est par exemple (possible mais) tr\`es
p\'enible de se convaincre que ce produit est, pour des alg\`ebres
commutatives, compatible (par passage aux quotients)
avec le produit usuel sur les complexes de De Rham commutatifs
(comparer avec notre remarque \ref{rem:DeRham}).

\section{Coextension de longueur 2 pour $sh$}\label{section:coext-long-2}

Le shuffle cyclique (\S \ref{subsection:shuffle}) donne une coextension de
longueur 2 du shuffle $\bbar{sh}$, c'est-\`a-dire que
l'application
 $$\bbar{sh}':\bar C_p(A)\otimes \bar C_q(B)\to \bar
C_{p+q+2}(A\otimes B)$$  (naturelle en les $k$-alg\`ebres $A$ et
$B$) v\'erifie~:
$$[\bbar{sh}, B]+[\bbar{sh}',b]=0,\qquad [\bbar{sh}',B]=0.$$
Nous verrons (\S \ref{section:genericite}) comment \'etendre ces \'equations \`a
la cat\'egorie des $k$-modules cycliques, et surtout (\S \ref{section:denormalisation}) comment
en d\'eduire, sur les complexes {\sl non normalis\'es\/}, une
coextension de $sh$, le prix \`a payer \'etant de perdre la
propri\'et\'e \og\ de longueur 2\fg\ (remarque \ref{rem:V.5}). Cependant, si l'anneau $k$ (est de
caract\'eristique $0$ et) contient ${\Bbb Q}$, cette longueur 2
est pr\'eserv\'ee par une construction directe. L'objet du pr\'esent
paragraphe est en effet la proposition suivante.

\begin{prop}\label{prop:coextsh}Sur les ${\Bbb Q}$-modules cycliques, $sh$ admet une coextension naturelle de longueur 2.
\end{prop}

\noindent{\sl Preuve.\/} La m\'ethode est un raffinement de celle
de Hood-Jones (\S \ref{subsection:Hood-Jones}). Le but est de construire des applications
$sh'_n:(P_*\otimes Q_*)_n\to(P\times Q)_{n+2}$ ($n\geq 0$)
naturelles en les ${\Bbb Q}$-modules cycliques $P,Q$, v\'erifiant~:
\begin{equation*}bsh'_n=Z_n:=sh'_{n-1}b+sh_{n+1}B-Bsh_n\qquad (n\geq 0)\tag*{$(1)_n$}\end{equation*}
\begin{equation*}sh'_nB=Bsh'_{n-1}\qquad(n\geq 1)\tag*{$(2)_n$}\end{equation*}
(avec $(P_*\otimes Q_*)_{-1}=0$, donc $sh'_{-1}=0$).

Posons
$$T_k:=1-(-1)^kt_k$$ et rappelons que
$$ B=TrN,\qquad{\rm avec}\qquad r_k:=t_{k+1}s_k\qquad{\rm et}\qquad
N_k=\sum_{i=0}^k(-1)^{ki}t_k^i.$$
J'en d\'eduis l'\'equation~:
$$d_0 B=2N.$$
(Si l'on pr\'ef\`ere les conventions de Connes, il faut choisir
$r_k=(-1)^ks_k$, et donc remplacer $d_0$ par $(-1)^kd_0t_{k+1}=b'-b$).
Sur ${\Bbb Q}$, on dispose classiquement d'op\'erateurs $h'_k:={1\over k+1}$ et
$h_k$ tels que $$hT=Th=1-h'N.$$ Posons alors
$$ L:=h'{d_0\over 2},\qquad K:=1-BL.$$
On remarque que $hT$, $h'N$, $K$, $1-K$ sont idempotents, et v\'erifient les
relations utiles suivantes~:
$$h'N=LB,\qquad KB=0,\qquad hT(1-K)=1-K, \qquad hT^2=T.$$

Les $sh'_n$ vont \^etre construits par r\'ecurrence sur $n$. Il est
n\'ecessaire pour cela d'ajouter l'\'equation
\begin{equation*}Bsh'_n(T\otimes T)=0\qquad (n\geq 2)\tag*{$(3)_n$}\end{equation*}
car $(2)_n\Rightarrow (3)_{n-1}$. (Pour $n=0$ ou $1$, $T\otimes T=0$ car
$T_0=0$).

Pour $n\geq 2$, voici l'\'etape de r\'ecurrence de la construction de $sh'_n$,
c'est-\`a-dire d'\'el\'ements $sh'_{p,q}\in
\lambda  [p]_{n+2}\otimes \lambda [q]_{n+2}$ pour $p+q=n$, v\'erifiant les
\'equations
\begin{equation*}bsh'_{p,q}=Z_{p,q}\qquad (0\leq p\leq
n)\tag*{$(1)_{p,q}$}\end{equation*}
\begin{equation*}sh'_{p,q}(B\otimes 1)=\alpha
:=(-1)^psh'_{p-1,q+1}(1\otimes B)+Bsh'_{p-1,q}\qquad(1\leq p\leq n)\tag*{$(2)_{p,q}$}\end{equation*}
\begin{equation*}Bsh'_{p,q}(T\otimes T)=0\qquad(1\leq p\leq n-1)\tag*{$(3)_{p,q}$}\end{equation*}
auxquelles il est n\'ecessaire de rajouter
\begin{equation*}sh'_{p,q}(T\otimes B)=\beta (T\otimes 1)\text{ avec }\beta
:=(-1)^pBsh'_{p,q-1}\quad (1\leq p\leq n-1)\tag*{$(4)_{p,q}$}\end{equation*}
car $(2)_{p,q}\Rightarrow (4)_{p-1,q+1}$.

Pour $n$ fix\'e, les $sh'_{p,q}$ ($p+q=n$) sont construits par r\'ecurrence
sur $p$. Par hypoth\`ese de r\'ecurrence sur $n$, $bZ_{p,q}=0$ donc (par
acyclicit\'e de $\Lambda [p]\times\Lambda [q]$ en degr\'e $n+1\geq 3$) il
existe $Y_{p,q}\in \lambda [p]_{n+2}\otimes \lambda [q]_{n+2}$
tel que $bY_{p,q}=Z_{p,q}$. On peut donc poser $sh'_{0,n}=Y_{0,n}$, mais pour
$p>0$ il faut modifier $Y_{p,q}$ de fa\c con \`a satisfaire aussi les
\'equations $(2)_{p,q}$ et (si $p<n$) $(3)_{p,q}$ et $(4)_{p,q}$. On va pour
cela chercher une solution $\gamma $ de ces trois \'equations, puis la
\og~recoller\fg\ avec $Y_{p,q}$.

Simplifions d'abord $(2)_{p,q}$ et $(4)_{p,q}$. Par hypoth\`ese de
r\'ecurrence, $\alpha (T\otimes 1)=0$ (d'apr\`es $(4)_{p-1,q+1}$) et $\beta
(T\otimes T)=0$ (d'apr\`es $(3)_{n-1}$). Donc $\alpha =\alpha _0(B\otimes 1)$
et $\beta (T\otimes 1)=\beta _0(T\otimes B)$, en posant $\alpha _0:=\alpha
(L\otimes 1)$ et $\beta _0=\beta (1\otimes L)$, d'o\`u~:
$$\begin{matrix}(2)_{p,q}&\Longleftrightarrow &(sh'_{p,q}-\alpha _0)(B\otimes
1)=0&{\rm\ et}\\
(4)_{p,q}&\Longleftrightarrow &(sh'_{p,q}-\beta _0)(T\otimes
B)=0.&\end{matrix}$$
Posons alors $$\gamma =\alpha _0((1-K)\otimes 1)+\beta _0(hTK\otimes(1-K)).$$
Ainsi, $\gamma $ est solution de $(2)_{p,q}$, mais aussi de $(3)_{p,q}$ (car par
construction,\break $B\alpha (1\otimes T)=0$ et $B\beta =0$). D'autre part
$\alpha (1\otimes B)=\beta (B\otimes 1)$ (d'apr\`es $(2)_{n-1}$), d'o\`u
$\alpha _0(B\otimes B)=\beta _0(B\otimes B)$, ou encore~:
$\alpha _0((1-K)\otimes(1-K))=\beta _0((1-K)\otimes(1-K))$, si bien que
$\gamma $ peut aussi s'\'ecrire~:
$$\gamma =\alpha _0((1-K)\otimes K)+\beta _0(hT\otimes(1-K)),$$
qui est solution de $(4)_{p,q}$.

\`A pr\'esent, \og~recollons\fg\  $\gamma $ avec $Y_{p,q}$. Par
hypoth\`ese de r\'ecurrence,\break
$BZ(T_p\otimes T_q)=0$ (d'apr\`es
$(3)_{n-1}$) d'o\`u, en posant $Y'=Y_{p,q}(Th\otimes Th)$~: $BbY'=0$, donc
$NbY'={d_0\over 2}BbY'=0$ donc (puisque $H_{n+2}^\lambda(\Lambda
[p]\otimes\Lambda [q])=0$) il existe $U,V$ tels que $Y'=TU+bV$. Posons
$$sh'_{p,q}:=\gamma +TU(ThK\otimes ThK)+Y(ThK\otimes h'N+h'N\otimes 1).$$
Les \'equations $(2)-(3)-(4)_{p,q}$, satisfaites par $\gamma $, le sont
encore par $sh'_{p,q}$. Reste \`a v\'erifier $(1)_{p,q}$. Par hypoth\`ese de
r\'ecurrence, $b\alpha =Z(B\otimes 1)$ et $b\beta (T\otimes 1)=Z(T\otimes B)$
(d'apr\`es $(1)_{n-1}$, $(2)_{n-1}$, et $(1)_{p-1,q+1}$). On en d\'eduit
facilement que $b\gamma =Z((1-K)\otimes K+Th\otimes (1-K))$. D'autre part, par
choix de $U$, $bTU(ThK\otimes ThK)=bY(ThK\otimes ThK)$. Donc $sh'_{p,q}$ est
bien solution de $(1)_{p,q}$, puisque
$$\displaylines{(1-K)\otimes K+Th\otimes (1-K)+ThK\otimes ThK+ThK\otimes
h'N+h'N\otimes 1=\cr
(1-K)\otimes K+Th\otimes (1-K)+ThK\otimes K+h'N\otimes 1=\cr
Th\otimes (1-K)+Th\otimes K+h'N\otimes 1=\cr
Th\otimes 1+h'N\otimes 1=1.\cr}$$

Pour finir, initialisons la r\'ecurrence sur $n$, en montrant qu'il existe
$sh'_{0,0}$, $sh'_{1,0}$, $sh'_{0,1}$ v\'erifiant les \'equations $(1)_{0,0}$,
$(1)_{0,1}$, $(1)_{1,0}$ et $(2)_{1,0}$. Le seul probl\`eme est l'\'equation
$(2)_{1,0}$ puisque les \'equations $(1)_n$ correspondent simplement \`a la
recherche d'une coextension de $sh$ de longueur quelconque (\'eventuellement
infinie), qui sera r\'esolue au \S \ref{section:denormalisation}. Autrement dit~: il nous reste \`a
modifier une solution donn\'ee $Y_0, Y_1$ (cf remarque V.5) de $(1)_0$ et
$(1)_1$ de fa\c con \`a obtenir $sh'_1B=Bsh'_0$. Il suffit pour cela de modifier
$Y_{1,0}$ comme dans le cas g\'en\'eral ci-dessus. En effet, dans cette \'etape
de la construction par r\'ecurrence, le fait que $n\geq 2$ n'avait pas servi.
Plus explicitement~: on pose $sh'_{0,0}:=Y_{0,0}$, $sh'_{0,1}:=Y_{0,1}$ et
$sh'_{1,0}:=\alpha _0((1-K)\otimes 1)+Y_{1,0}(K\otimes 1)$ avec $\alpha
_0=\alpha (L\otimes 1) $ et $\alpha =-Y_{0,1}(1\otimes B)+BY_{0,0}$. En effet,
dans le cas $(p,q)=(1,0)$, la solution g\'en\'erale est simplifi\'ee par les
remarques suivantes~: $\beta =0$ (donc $\beta _0=0$), $T_1\otimes T_0=0$ (donc
$Y'=0$, donc $U=0, V=0$ conviennent), $h'_0N_0=1$ (donc $ThK\otimes
h'_0N_0+h'N\otimes 1=K\otimes 1$) et $K_0=1$.\cqfd

\section{G\'en\'ericit\'e des modules cycliques
d'alg\`ebres}\label{section:genericite}

Les coextensions $\bbar{sh}_\infty, \bbar{AW}_\infty$ du \S \ref{section:coext} seront
construites \`a partir d'une homotopie naturelle $\varphi $ entre $\bbar{sh}\
\bbar{AW}$ et ${\rm id}$. Pour trouver $\varphi$, nous utiliserons les
\og~mod\`eles acycliques~\fg\ de Kassel (\S \ref{subsection:Kassel}) plut\^ot que ceux de
Hood-Jones (\S \ref{subsection:Hood-Jones}), donc nous nous restreindrons aux modules
cycliques d'alg\`ebres (ceci afin de disposer de l'outil des
r\'esolutions pour construire les homotopies).

Mais le $\varphi:\bar C_n(A\otimes B)\to \bar C_{n+1}(A\otimes B)$ que nous
obtiendrons (naturel en les alg\`ebres $A$ et $B$) sera
\og~cyclique~\fg, c'est-\`a-dire exprimable (lin\'eairement) en termes de
morphismes de la cat\'egorie $\Lambda ^{op}$. Or l'objet du
pr\'esent paragraphe peut s'\'enoncer informellement comme suit.

\begin{princ}\label{princ}Si une propri\'et\'e
cyclique est vraie pour les modules cycliques d'alg\`ebres alors
elle est vraie pour tous les modules cycliques.
\end{princ}

\noindent\hfill En particulier, nous disposerons d'une homotopie $\varphi
:\bbar{(P\times Q)}_n\to\bbar{(P\times Q)}_{n+1}$\break (naturelle
en les modules cycliques $P,Q$) entre $\bbar{sh}\ \bbar{AW}$ et ${\rm
id}$, comme r\'esultat du \S III.1 ci-dessous. Le \S III.2 ne sera
pas utilis\'e dans la suite, mais illustre une autre facette du
principe III.1.

\subsection{Prolongement des \'egalit\'es}

Ce principe sera pr\'esent\'e
dans le cas de $k$-modules simpliciaux, pour all\'eger l'expos\'e,
mais la preuve s'adapte sans aucune difficult\'e aux $k$-modules
cycliques, et aux produits de deux tels modules. Soit $F:C_m(A)\to
C_n(A)$ une application naturelle en la $k$-alg\`ebre $A$. Notre
but est d'identifier les $F$ qui s'\'etendent en une application
$f:P_m\to P_n$ naturelle en le $k$-module simplicial $P$ (un tel
$F$ sera dit {\sl $\Delta $-repr\'esentable\/}), de montrer
qu'alors le repr\'esentant $f$ est unique, et de \og~lire sur
$F$~\fg\ \`a quelle condition $f$ passe aux quotients,
c'est-\`a-dire d\'efinit une application $\bar f:\bar P_m\to \bar
P_n$ ($f$ sera alors dite {\sl\og~normalisable~\fg\/}) . On
identifiera $f$ avec son repr\'esentant canonique $f(1_m)\in
k[\Delta ^{op}(m,n)]=k[\Delta [m]_n]$, et de m\^eme $F$ avec
$F(a_0\otimes\ldots\otimes a_m)\in C_n(T_m)$.

Remarquons que $C_n(T_m)=k[M_m^{n+1}]$, o\`u $M_m$ est le mono\"\i
de libre sur l'en\-semble $\{a_0,\ldots,a_m\}$~; ou encore~:
$$C_n(T_m)=k[C_n(M_m)],$$
en appelant $C_*(M)$ l'ensemble simplicial naturellement associ\'e
\`a un mono\"\i de $M$. Etudier le morphisme de $k$-modules
(libres) $$\begin{matrix}
k[\Delta [m]_n]&\to&C_n(T_m)\\
f&\mapsto&F=f(a_0\otimes\ldots\otimes a_m)\end{matrix}$$
 se ram\`ene donc \`a \'etudier l'application
 $$\Delta [m]_n\to C_n(M_m)$$
 (c'est-\`a-dire en fait~: le morphisme d'ensembles simpliciaux
$\Delta [m]\to C(M_m)$ canoniquement associ\'e, par propri\'et\'e universelle
de $\Delta [m]$, \`a l'\'el\'ement $(a_0,\ldots,a_m)$ de $C_m(M_m)$).

\begin{lem}
\begin{enumerate}[(i)]
\item L'application naturelle $\Delta [m]_n\to C_n(M_m)$ est
injective.
\item Son image est constitu\'ee des \'el\'ements de la forme
$$(m_0,\ldots,m_n)=(a_0\ldots a_{i_0}, a_{i_0+1}\ldots
a_{i_1},\ldots , a_{i_{n-1}+1}\ldots a_{i_n})$$
avec $0\leq i_0\leq i_1\leq\ldots\leq i_n=m$.
\item Un tel \'el\'ement est l'image d'un \'el\'ement
d\'eg\'en\'er\'e de $\Delta [m]_n$ si et seulement s'il existe
$\ell>0$ tel que $m_\ell=1$.
\end{enumerate}
\end{lem}

\noindent{\sl Preuve.\/} Un \'el\'ement de $\Delta [m]_n$
s'\'ecrit de mani\`ere unique sous la forme $$f=s_{k_1}\ldots
s_{k_p}d_{n-p}^{j_{n-p}}\ldots d_0^{j_0}$$
avec $p\geq 0$ ($p>0$ si et seulement si $f$ est
d\'eg\'en\'er\'e), $0=:k_0\leq k_1\leq\ldots\leq k_{p+1}:=n-p$,
$j_0,\ldots, j_{n-p}\geq 0$, et $n+j_0+\ldots+j_{n-p}=m+p$.
L'image de $f$ dans $C_n(M_m)$ est alors de la forme \'enonc\'ee
dans {\sl i)\/} avec, pour $k_t+t\leq s\leq k_{t+1}+t$,
$i_s=j_0+\ldots+j_{s-t}+s-t$. En particulier, les $\ell$ tels que
$m_\ell=1$ sont les $\ell$ de la forme $k_t+t$ pour $1\leq t\leq
p$, d'o\`u {\sl iii)\/}. R\'eciproquement un
$(m_0,\ldots,m_n)$ comme dans {\sl i)\/} admet un unique
ant\'ec\'edent $f$ (d'o\`u {\sl i)\/} et {\sl ii)\/})~: d'abord,
$p$ et $(k_1,\ldots, k_p)$ son d\'etermin\'es par la suite
$\ell_1<\ldots<\ell_p\leq n$ des indices $\ell>0$ tels que
$m_\ell=1$~; on peut ensuite se ramener au cas $p=0$, or dans ce
cas, $j_0,\ldots,j_n$ sont simplement d\'etermin\'es par~:
$j_0=i_0$ et pour $1\leq s\leq n$, $j_s=i_s-i_{s-1}-1$.\cqfd

\begin{prop}
\begin{enumerate}[(i)]
\item Un \'el\'ement $F\in C_n(T_m)$
est $\Delta $-repr\'esentable si et seulement s'il est combinaison
lin\'eaire d'\'el\'ements de la forme
$$(m_0,\ldots,m_n)=(a_0\ldots a_{i_0}, a_{i_0+1}\ldots
a_{i_1},\ldots , a_{i_{n-1}+1}\ldots a_{i_n})$$
avec $0\leq i_0\leq i_1\leq\ldots\leq i_n=m$, et
\item son
repr\'esentant $f$ est alors unique.
\item De plus, $f$ est normalisable
si et seulement si
$$\forall i\in\{1,\ldots,m\},\qquad \bbar{F(a_0\otimes\ldots\otimes
a_{i-1}\otimes 1\otimes a_{i+1}\otimes\ldots\otimes a_m)}=0$$
dans $\bar C_n(T_m)$, donc si et seulement si $F$ est
normalisable,
\item et dans ce cas, $\bar f=0$ si et seulement si $\bar
F=0$.
\end{enumerate}
\end{prop}

\noindent{\sl Preuve.\/} Les points {\sl i)\/} et {\sl ii)\/}
r\'esutent imm\'ediatement des points correspondants du
lemme. Dans les point {\sl iii)\/} et {\sl iv)} la partie
\og~seulement si~\fg\ est imm\'ediate. Reste donc \`a prouver~:

\begin{itemize}
\item{iv')} si $F\in k[X]$, o\`u $X$ d\'esigne l'ensemble des
\'el\'ements d\'eg\'en\'er\'es de $C_n(M_m)$, alors $f$ est de la
forme $\sum s_j g_j$ avec $g_j\in k[\Delta ^{op}(m,n-1)]$
\item{iii')} m\^eme chose en rempla\c cant $F$ (resp. $f$) par les
$F\circ s_i$ (resp. $f\circ s_i$) et $m$ par $m-1$.
\end{itemize}

Il suffit \'evidemment de prouver {\sl iv')\/}.
Supposons donc $F\in k[X]$, et notons $Y$ l'image de $\Delta [m]_n$
dans $C_n(M_m)$. Alors $F\in k[X]\cap k[Y]=k[X\cap Y]$. Or
d'apr\`es le point {\sl iii)\/} du lemme, $X\cap Y$ est l'image
dans $C_n(M_m)$ des \'el\'ements d\'eg\'en\'er\'es de $\Delta
[m]_n$, c'est-\`a-dire des \'el\'ements de la forme $s_j\circ g$
avec $g\in\Delta [m]_{n-1}$. Il existe donc des $g_j\in
k[\Delta ^{op}(m,n-1)]$ tels que $F=\sum s_j\circ
g_j(a_0\otimes \ldots\otimes a_m)$ donc (d'apr\`es le point
{\sl ii)\/} de la proposition) tels que $f=\sum s_j\circ
g_j$.\cqfd

R\'esumons la proposition pr\'ec\'edente~:

\begin{coro}\label{coro:reprF}
Soit $F:C_m(A)\to C_n(A)$
une application naturelle en l'alg\`ebre $A$ telle que
$F(a_0\otimes\ldots\otimes a_m)$ soit combinaison lin\'eaire
(\`a coefficients fix\'es) d'\'el\'ements de la forme
$$a_0\ldots a_{i_0}\otimes a_{i_0+1}\ldots
a_{i_1}\otimes\ldots\otimes a_{i_{n-1}+1}\ldots a_m.$$ Alors $F$ admet un
$\Delta $-repr\'esentant (unique) $f:P_m\to P_n$ (naturel en le module
simplicial $P_*$). Si de plus $F$ est normalisable alors $f$ l'est aussi, et
$\bar F=0\Longrightarrow\bar f=0$.
\end{coro}

On laisse au lecteur l'exercice de d\'emontrer, par la m\^eme
m\'ethode, les trois corollaires suivants.

\begin{coro}
Soit $F:C_m(A)\to C_n(A)$ une
application naturelle en l'alg\`ebre $A$ telle que
$F(a_0\otimes\ldots\otimes a_m)$ soit combinaison lin\'eaire
d'\'el\'ements de la forme
$$a_{i_{n-q}+1}\ldots a_{i_{n-q+1}}\otimes\ldots\otimes\ldots
a_m\otimes a_0\ldots a_{i_0}\otimes a_{i_0+1}\ldots
a_{i_1}\otimes\ldots\otimes\ldots a_{i_{n-q}}.$$
Alors $F$ admet un
$\Lambda $-repr\'esentant (unique) $f:P_m\to P_n$ (naturel en le module
cyclique $P_*$). Si de plus $F$ est normalisable
alors $f$ l'est aussi, et $\bar F=0\Longrightarrow\bar f=0$.
\end{coro}

\begin{coro}\label{coro:III.6}
Soit $F:C_m(A\otimes B)\to
C_n(A\otimes B)$ une application naturelle en les alg\`ebres $A,B$ telle que
$F((a_0\otimes b_0)\otimes\ldots\otimes (a_m\otimes b_m))$ soit
combinaison lin\'eaire d'\'el\'ements de la forme
$$\displaylines{(a_{i_{n-q}+1}\ldots
a_{i_{n-q+1}}\otimes\ldots\otimes\ldots a_m\otimes a_0\ldots
a_{i_0}\otimes a_{i_0+1}\ldots a_{i_1}\otimes\ldots\otimes\ldots
a_{i_{n-q}})\otimes\hfill\cr \hfill(b_{j_{n-r}+1}\ldots
b_{j_{n-r+1}}\otimes\ldots\otimes\ldots b_m\otimes b_0\ldots
b_{j_0}\otimes b_{j_0+1}\ldots b_{j_1}\otimes\ldots\otimes\ldots
b_{j_{n-r}}).\cr}$$ Alors $ F$ admet un
$\Lambda $-repr\'esentant (unique) $f:(P\times Q)_m\to
(P\times Q)_n$ (naturel en les modules cycliques
$P,Q$). Si de plus $F$ est normalisable
alors $f$ l'est aussi, et $\bar F=0\Longrightarrow\bar f=0$.
\end{coro}

\begin{coro}
Soit $F:(C_*(A)\otimes C_*(B))_m\to
(C_*(A)\otimes C_*(B))_n$ une application naturelle en les alg\`ebres $A,B$
admettant un $\Lambda $-repr\'esentant (n\'ecessairement unique) $f:(P\times
Q)_m\to (P\times Q)_n$ (naturel en les modules cycliques
$P,Q$). Si $F$ est normalisable
alors $f$ l'est aussi, et $\bar F=0\Longrightarrow\bar f=0$.
\end{coro}

\begin{rem}\label{rem:normalisable}
Une simple observation comme dans
\cite{G} (proof of theorem II.4.2) montre que dans les corollaires 4 et
5, $F$ (ou $f$) est automatiquement
normalisable d\`es que $n>m$, et $\bar F$ (ou $\bar f$) est automatiquement
nulle d\`es que $n>m+1$. Dans le corollaire 6, ces conditions sont \`a
remplacer respectivement par $n>2m$ et $n>2m+2$, et dans le corollaire 7, par
$n>m+1$ et $n>m+2$. On trouverait d'ailleurs les m\^emes conditions en
rempla\c cant 6 par \og~6bis~\fg\ avec $F:C_m(A\otimes B)\to
(C_*(A)\otimes C_*(B))_n$, et 7 par \og~7bis~\fg\ avec $F:(C_*(A)\otimes C_*(B))_m\to
C_n(A\otimes B)$. En particulier toute coextension de $sh$ est
normalisable, en une coextension de $\bbar{sh}$, et toute coextension de $\bbar
{sh}$ est de longueur 2.
\end{rem}

\subsection{G\'en\'eralisation de certains r\'esultats
de Cuntz-Quillen}

Les corollaires 6 et 7 ci-dessus s'appliqueront respectivement dans les
lemmes \ref{lem:barshbarAW} et \ref{lem:barAWBbarsh}~; le corollaire 5 va s'appliquer ici. Le
pr\'esent paragraphe est une digression qui illustre une autre facette du
principe \ref{princ}, en \'etendant aux modules cycliques des r\'esultats de
\cite{CQ1}, \cite{CQ2} sur les modules cycliques d'alg\`ebres.

Rappelons (preuve de la proposition \ref{prop:coextsh}) la d\'efinition de $r$ ($\bar r$ d\'esignera sur
$\Omega(A)$ l'op\'erateur usuellement not\'e ${\rm d}$)~:
$$r_n=t_{n+1}s_n,$$ et d\'efinissons l'op\'erateur de Karoubi
$\kappa $ par~: $$br+rb=1-\kappa $$ (donc $\kappa $ commute \`a
$b$).

\begin{prop}
$r$ et $\kappa $ sont
normalisables, $\bar r^2=0$, $[\bar\kappa ,\bar r]=0$,
$$\begin{matrix}\bar\kappa_{n+1}\bar r_n=(-1)^n\bbar{r_nt_n},&
\bar\kappa_{n+1}^{n+1}\bar r_n=\bar r_n,&
\bar\kappa _n^n=1+ b_{n+1}\bar \kappa _{n+1}^n\bar r_n,\\
\bar\kappa _{n-1}^nb_n=b_n,&
\bar \kappa _n^{n+1}=1-\bar r_{n-1}b_n,&
(\bar\kappa_n-1)(\bar\kappa _n^{n+1}-1)=0,\end{matrix}$$
$$\begin{matrix}
B_n=\sum_{j=0}^n\bar\kappa _{n+1}^j\bar r_n,&
\bar rB=B\bar r=B^2=0,\\
B\bar \kappa =\bar\kappa B=B,&
\bar\kappa_n^{n(n+1)}-1=bB=-Bb.\end{matrix}$$
\end{prop}

\noindent{\sl Preuve.\/} Ces propri\'et\'es ont
\'et\'e d\'emontr\'ees (\cite{CQ1} p.~81--83 ou \cite{CQ2} p.~387--389) pour
$P_*=C_*(A)$ ($A$ alg\`ebre unitaire) donc se g\'en\'eralisent \`a
tout module cyclique $P_*$ d'apr\`es le corollaire \ref{coro:reprF}. On peut
bien s\^ur aussi les red\'emontrer par un calcul direct dans
$\Lambda ^{op}$. Par exemple~: $r$ est normalisable car $\forall
i\in\{0,\ldots, n\}$, $r_{n+1}s_i=s_{i+1}r_n$~; $\bar r^2=0$ car
$r_{n+1}r_n=s_0r_n$~; $\kappa $ est normalisable parce que $r$ et
$b$ le sont (ou directement~: parce que
$$\kappa _n=(-1)^nt_n(1-s_{n-1}d_n),$$
d'o\`u $\kappa _{n+1}s_i=-s_{i+1}\kappa _n$ pour $i<n$ et
$\kappa _{n+1}s_n=0$). Mais le calcul dans $\Lambda ^{op}$ devient
plus p\'enible au fil des \'equations, d'o\`u l'int\'er\^et du
principe \ref{princ}.\cqfd

Si $k$ contient ${\Bbb Q}$, on
peut encha\^\i ner sur la \og~d\'e\-com\-po\-si\-tion spectrale
relative \`a $\bar \kappa $~\fg\ en rempla\c cant partout, dans \cite{CQ1}
(p.~84-85) et \cite{CQ2} (p.~389-392), $\Omega(A)$ par $\bar
P_*$, pour un module cyclique quelconque $P_*$. Mais une propri\'et\'e
essentielle de $\Omega(A)$, utilis\'ee dans les cons\'equences de cette
d\'ecomposition, est~: $H_n(\Omega(A),\bar r)=0$ pour $n\geq 1$. Or {\sl cette
propri\'et\'e ne se g\'en\'eralise pas \`a tout module cyclique,\/} ce qui
met en \'evidence (sans le contredire) les limites du principe \ref{princ}. Les
modules cycliques $\lambda [n]$ ne fournissent pas un contre-exemple \`a
cette propri\'et\'e, donc la m\'ethode la plus naturelle, pour
construire un module cyclique $P$ tel que $H_1(\bar P,\bar r)\neq 0$, est de
quotienter $\lambda [1]$ en \og~for\c cant~\fg\ $r_1(1_1)$ \`a \^etre
d\'eg\'en\'er\'e, et de v\'erifier que dans ce quotient $P$, $[
1_1]\notin{\rm Im}(\bar r_0)$.

\begin{prop}
Soit $P$ le module cyclique quotient de $\lambda [1]$ par la
relation~: $(1-s_1d_1)(1-s_0d_0)r_1(1_1)=0$. Alors $H_1(\bar P,\bar r)\neq 0$.
\end{prop}

\noindent{\sl Preuve.\/} Nous allons montrer bien plus que la proposition, en
explicitant com\-pl\`e\-te\-ment $P$. Il faut pour cela d\'ecrire le
sous-module cyclique $Q$ de $\lambda [1]$ engendr\'e par l'\'el\'ement
$(1-s_1d_1)(1-s_0d_0)r_1(1_1)$. Introduisons quelques notations~:
$F:=d_0(1_1)$, $G:=d_1(1_1)$, $E_i:=s_0^i(1_1)$, $s:=s_n\in\Lambda
^{op}(n,n+1)$,
$$\begin{matrix}f_1:=&sF+tsG-(1+t)E_0;\\
\text{si }n\geq 2, f_n:=&s^nF+t^n(s^{n-2}E_1-s^{n-1}E_0)-s^{n-1}E_0,\\
e_{n,n}:=
&s^nG-(1+t+\ldots+t^{n-1})(s^{n-2}E_1-s^{n-1}E_0)-s^{n-1}E_0;\\
\text{si }0<i<n, e_{i,n}:=
&s^{n-1-i}E_i-(1+\ldots+t^{i-1})(s^{n-2}E_1-s^{n-1}E_0)-s^{n-1}E_0.\end{matrix}$$
(En particulier $e_{1,n}=0$). Soit $R$ le sous-module de $\lambda [1]$
lin\'eairement engendr\'e par les $f_n$, $e_{i,n}$ et leurs images par les
puissances de $t$. Nous allons prouver que $Q=R$. On a d\'ej\`a
$(1-s_1d_1)(1-s_0d_0)r_1(1_1)=-tf_2\in R$. Il suffit ensuite de v\'erifier
qu'\`a l'aide des $s_j$, $d_j$, $t$, on peut engendrer, \`a partir de $f_2$,
tous les $f_n$, $e_{i,n}$ (ce qui donnera $R\subset Q$), et que
toute image par $s_j$ ou $d_j$ d'un
$f_n$ ou d'un $e_{i,n}$ appartient \`a $R$ (ce qui prouvera que $R$ est un
sous-module cyclique).\\
Ces propri\'et\'es se d\'eduisent facilement des
\'equations suivantes :
$$d_0f_n=d_nf_n=0,\quad s_0f_n=f_{n+1}+t^{n+1}e_{2,n+1},\quad s_nf_n=f_{n+1}+t^ne_{2,n+1},$$
$$0<j<n\Rightarrow\left[d_jf_n=f_{n-1}\text{ et }s_jf_n=f_{n+1}\right]$$
et pour $2\leq i\leq n$ :
$$s_0(e_{i,n})=e_{i+1,n+1}-e_{2,n+1},\quad d_0(e_{i,n})=e_{i-1,n-1},$$
$$0<j<i\Rightarrow\left[s_j(e_{i,n})=e_{i+1,n+1}-(t^{j-1}+t^j)e_{2,n+1}\text{ et }d_j(e_{i,n})=e_{i-1,n-1}+t^jf_{n-1}\right],$$
$$s_i(e_{i,n})=e_{i+1,n+1}-t^{i-1}e_{2,n+1},\quad d_i(e_{i,n})=e_{i-1,n-1},$$
$$j>i\Rightarrow\left[s_j(e_{i,n})=e_{i+1,n+1}\text{ et }d_j(e_{i,n})=e_{i,n-1}\right].$$
On a donc $P=\lambda [1]/R$. Or tout
\'el\'ement de $\Lambda [1]$ s'\'ecrit de fa\c con unique sous la forme
$t^ms^kE_i$ ou $t^ms^kF$ ou $t^ms^kG$. Vue la forme des g\'en\'erateurs de
$R$, on en d\'eduit donc une {\sl base\/} de $P$~:
$$\displaylines{P_0=k[F]\oplus k[G],\qquad
P_1=k[E_0]\oplus k[tE_0]\oplus k[sF]\oplus k[sG],{\rm \ et }\cr
{\rm pour\  }n\geq 2,\qquad
P_n=\bigoplus_{0\leq m\leq n}(k[t^ms^{n-1}E_0]\oplus k[t^ms^{n-2}E_1]).\cr}$$
Par construction, $\bar r_1([E_0])=0$, mais $[E_0]$ n'est pas combinaison
lin\'eaire de $r_0[F]=[E_0]+[tE_0]-[sG]$, de $r_0[G]=[E_0]+[tE_0]-[sF]$, et
des d\'eg\'en\'er\'es $[sF]$, $[sG]$.\cqfd

\section{Coextensions explicites de $\bbar{sh}$ et
$\bbar{AW}$.}\label{section:coext}

Les ($b$-) quasi-isomorphismes $\bbar{sh}$ et $\bbar{AW}$
(inverses l'un de l'autre \`a homotopie pr\`es) v\'erifient les
hypoth\`eses du lemme de perturbation (lemme \ref{lem:perturbation}) car on a de
plus~:
$$\bbar{AW}\ \bbar{sh}=1.$$
Malheureusement, aucun des deux ne commute \`a $B$ (sinon, le
probl\`eme de coextension serait imm\'ediatement r\'esolu, par le
m\^eme raisonnement qu'au \S \ref{subsection:perturbnorma}). Mais nous sommes sauv\'es par
une propri\'et\'e providentielle (et \`a ma connaissance,
ignor\'ee jusqu'\`a pr\'esent)~:

\begin{lem}\label{lem:barAWBbarsh}
On a $\bbar{AW}B\bbar{sh}=B$.
\end{lem}

\noindent{\sl Preuve.\/} Soient $x=a_0{\rm d}a_1\ldots{\rm
d}a_p\in\Omega^p(A), y=b_0{\rm d}b_1\ldots{\rm
d}b_q\in\Omega^q(B)$ (pour deux alg\`ebres unitaires quelconques
$A,B$). Il suffit (corollaire III.7) de prouver que
$\bbar{AW}B\bbar{sh}(x\otimes y)=B(x\otimes y)$.
Or $\bbar{sh}(x\otimes y)=\sum_{\sigma \in S_{p,q}}\varepsilon
(\sigma ) (a_0\otimes b_0)\sigma (z)$, o\`u $z$ d\'esigne
la suite ${\rm d}a_1,\ldots,{\rm
d}a_p,{\rm d}b_1,\ldots,{\rm d}b_q$ \`a laquelle s'applique le
$(p,q)$-shuffle $\sigma\in S_{p,q}$, donc $B\bbar{sh}(x\otimes
y)=\sum_{k=0}^{p+q}(-1)^{(p+q)k}\sum_{\sigma \in S_{p,q}}
\varepsilon (\sigma )t^k({\rm d}(a_0\otimes b_0),\sigma (z))$,
o\`u $t$ d\'esigne la permutation circulaire
$t(w_0,\ldots,w_{p+q})=(w_{p+q},w_0,\ldots,w_{p+q-1})$.\\
Puis on applique $\bbar {AW}_{p+q+1}$, qui va \og~tuer~\fg, dans
cette somme, tous les permut\'es o\`u un ${\rm d}b$ appara\^\i t
devant un ${\rm d}a$. En effet, $\bbar{AW}_n({\rm d}(u_1\otimes
v_1)\ldots{\rm d}(u_n\otimes v_n))=\sum_{k=0}^n(u_{k+1}\ldots
u_n{\rm d}u_1\ldots{\rm d}u_k)\otimes (v_1\ldots v_k{\rm
d}v_{k+1}\ldots{\rm d}v_n)$, donc s'il existe un $i<j$ tel que
$u_i=1$ et $v_j=1$ alors les termes pour $k\geq i$ seront nuls
parce que ${\rm d}u_i=0$, et ceux pour $k\leq i$ aussi parce que
${\rm d}v_j=0$.\\
Dans $\bbar{AW}B\bbar{sh}(x\otimes y)$,
$\bbar{AW}$ ne s'applique donc de fa\c con non nulle qu'\`a des produits de la
forme ${\rm d}a_{i_1}\ldots{\rm d}(a_0\otimes b_0)\ldots{\rm
d}a_{i_p}{\rm d}b_1\ldots{\rm d}b_{j_q}$ (dans ce produit, ${\rm
d}(a_0\otimes b_0)$ peut \^etre ins\'er\'e n'importe o\`u). Un tel
produit, figurant dans $B\bbar{sh}(x\otimes y)$, ne peut \^etre
qu'un permut\'e circulaire $t^k$ de
$$\begin{matrix}&{\rm d}(a_0\otimes b_0)\sigma _i(z)=&
{\rm d}(a_0\otimes b_0){\rm d}a_1\ldots{\rm d}a_i{\rm
d}b_1\ldots{\rm d}b_q{\rm d}a_{i+1}\ldots{\rm d}a_p\\
&{\rm avec\ }i\geq 1,k=p-i,\cr
{\rm ou\ de}&{\rm d}(a_0\otimes b_0)\tau
_j(z)=&{\rm d}(a_0\otimes b_0){\rm d}b_1\ldots{\rm d}b_j{\rm
d}a_1\ldots{\rm d}a_p{\rm d}b_{j+1}\ldots{\rm d}b_q\\
&{\rm avec\ }j\geq 1,k=p+q-j.&\end{matrix}$$
Donc $\bbar{AW}B\bbar{sh}(x\otimes y)=$ $$
\displaylines{\sum_{i=1}^p\varepsilon (\sigma
_i)(-1)^{(p+q)(p-i)}{\rm d}a_{i+1}\ldots{\rm d}a_p{\rm
d}a_0\ldots{\rm d}a_i\otimes b_0{\rm d}b_1\ldots{\rm
d}b_q\ +\hfill\cr\hfill\sum_{j=1}^q\varepsilon (\tau
_j)(-1)^{(p+q)(p+q-j)} a_0{\rm d}a_1\ldots{\rm d}a_p\otimes{\rm
d}b_{j+1}\ldots{\rm d}b_q{\rm d}b_0\ldots{\rm d}b_j.\cr}$$
$$\begin{matrix}\text{Or}&\varepsilon (\sigma
_i)(-1)^{(p+q)(p-i)}&=&(-1)^{q(p-i)+(p+q)(p-i)}&=&(-1)^{p(i+1)}\\
\text{et}&\varepsilon (\tau
_j)(-1)^{(p+q)(p+q-j)}&=&(-1)^{pj+(p+q)(p+q-j)}&=&(-1)^{q(j+1)+p},\end{matrix}$$
do\`u le r\'esultat~:
$$\bbar{AW}B\bbar{sh}(x\otimes y)=(Bx)\otimes y+(-1)^px\otimes
(By)=B(x\otimes y).\qquad\Box$$

\begin{lem}
Soit $\varphi $ une
homotopie (naturelle) telle que $\bbar{sh}\ \bbar{AW}=1+b\varphi
+\varphi b$. En appliquant le lemme \ref{lem:perturbation}, on obtient~:
$b_\infty=b+B$.
\end{lem}

\noindent{\sl Preuve.\/} La composante de degr\'e $-1$ de
$b_\infty$ est $b$ par d\'efinition. Celle de degr\'e $1$ est $B$
d'apr\`es le lemme \ref{lem:barAWBbarsh}. Les suivantes sont nulles quel que soit
le choix de $\varphi $, pour des raisons de degr\'e (remarque \ref{rem:normalisable}).\cqfd

Si de plus $\varphi $ est choisie sp\'eciale, le lemme \ref{lem:perturbation} donne
donc~:

\begin{thm}\label{thm:perturbcoext}
Le lemme de perturbation
fournit des coextensions $\bbar{AW}_\infty$ et $\bbar{sh}_\infty$,
mutuellement quasi-inverses en homologie cyclique n\'egative.
\end{thm}

\noindent{\sl Remarque.\/} On a d\'ej\`a vu (remarque \ref{rem:normalisable}) que
$\bbar{sh}_\infty$ \'etait automatiquement de longueur 2. Par
contre, $$\varphi _k:\Omega^n(A\otimes
B)\to\Omega^{n+2k+1}(A\otimes B)$$ n'est nul que pour $2k>n+1$, et
$$\bbar{AW}_k:\Omega^n(A\otimes B)\to
(\Omega(A)\otimes\Omega(B))^{n+2k}$$
n'est nul que pour $2k>n+2$.

Reste \`a construire un $\varphi $ de mani\`ere \`a rendre ce
th\'eor\`eme explicite, et en particulier \`a identifier
$\bbar{sh}_1=\varphi B\bbar{sh}$. Nous allons pour cela utiliser
une m\'ethode de mod\`eles acycliques comme dans \cite{HJ} (cf. notre \S
\ref{subsection:Hood-Jones}, et notre \S II qui s'en inspire), ou plut\^ot comme dans
\cite{K2} (cf. notre \S \ref{subsection:Kassel}), puisque nos \og~mod\`eles~\fg\ seront les
alg\`ebres libres $T_n$ (d\'ej\`a rencontr\'ees au \S \ref{section:genericite}).

On veut construire, par r\'ecurrence sur $n$,
$$\varphi _n:\Omega^n(A\otimes B)\to\Omega^{n+1}(A\otimes B)$$
naturelle en les alg\`ebres $A$ et $B$ -- donc repr\'esent\'ee par un
$\Phi _n\in\Omega^{n+1}(T_n\otimes T_n)$ -- telle
que $b\varphi _n=\bbar{sh}\ \bbar{AW}-1-\varphi _{n-1}b$
(avec $\varphi _{-1}=0$), ou encore~:\\
\centerline{$b(\Phi _n)=Z_n:=(\bbar{sh}\ \bbar{AW}-1-\varphi
_{n-1}b)(a_0\otimes b_0{\rm d}(a_1\otimes b_1)\ldots {\rm
d}(a_n\otimes b_n)).$}
Supposons $\varphi _0,\ldots,\varphi _{n-1}$ construits, avec $n\geq
2$ (le probl\`eme de l'initialisation sera r\'egl\'e facilement).
Pour construire $\Phi _n\in\Omega^{n+1}(T_n\otimes
T_n)$ il suffit de remarquer que par hypoth\`ese de
r\'ecurrence $b(Z_n)=0$, et d'utiliser l'acyclicit\'e de
$\Omega(T_n\otimes T_n)$ en degr\'es $\geq 3$. Mais {\sl pour v\'erifier que
le $\varphi _n$ associ\'e est normalisable, il faut l'expliciter\/} (sinon,
on retombe sur le probl\`eme \'evoqu\'e aux \S \ref{subsection:Hood-Jones} et \ref{subsection:Kassel}). Il faut donc d'abord
construire une homotopie $h$ qui t\'emoigne de cette acyclicit\'e de
$\Omega(T_n\otimes T_n)$.

Nous allons calculer $h$ \`a l'aide de r\'esolutions.
C'est ce qui nous a conduits \`a nous restreindre aux modules
cycliques d'alg\`ebres. Mais gr\^ace au corollaire \ref{coro:III.6},
le $\varphi $ que nous obtiendrons sera (vue sa forme) $\Lambda
$-repr\'esentable et son repr\'esentant v\'erifiera la m\^eme
\'equation, et de m\^eme pour les applications qui s'en
d\'eduisent dans le th\'eor\`eme IV.3.

\begin{lem}\label{lem:homotopie}
Soient $V, W$ deux $k$-modules et $A:=T(V)$,
$B:=T(W)$ leurs alg\`ebres tensorielles. On a $$hb+bh=1-k,$$ avec
$h_n:\Omega^n(A\otimes B)\to \Omega^{n+1}(A\otimes B)$, $k_n:\Omega^n(A\otimes
B)\to \Omega^n(A\otimes B)$ d\'efinies par~: $k_0:=1$, $k_n:=0$ pour
$n\geq 2$, et si $u_i, v_j\in V$, $w_k\in
W$, $a\in A$, $b\in B$ et $c\in A\otimes B$,
$$\displaylines{k_1(c {\rm d}(v_1\ldots v_q\otimes w_1\ldots w_r)):=
\sum_{k=1}^r (v_1\ldots v_q\otimes w_{k+1}\ldots
w_r)c(1\otimes w_1\ldots w_{k-1}){\rm d}(1\otimes w_k)\hfill\cr
\hfill+\sum_{j=1}^q (v_{j+1}\ldots v_q\otimes 1)c(v_1\ldots
v_{j-1}\otimes w_1\ldots w_r){\rm d}(v_j\otimes 1),\cr
k_2(c{\rm d}(v_1\ldots v_q\otimes b){\rm d}(a\otimes w_1\ldots
w_r)):=
\sum_{k=1}^r\sum_{j=1}^q (v_{j+1}\ldots v_qa\otimes
w_{k+1}\ldots w_r)c\hfill\cr\hfill(v_1\ldots v_{j-1}\otimes
bw_1\ldots w_{k-1}) [{\rm d}(v_j\otimes 1){\rm d}(1\otimes w_k)-
{\rm d}(1\otimes w_k){\rm d}(v_j\otimes 1)],\cr
h_n(\omega_{n-2}{\rm
d}(u_1\ldots u_p\otimes b){\rm d}(v_1\ldots v_q\otimes w_1\ldots
w_r)):=\hfill\cr (-1)^{n+1}\sum_{k=1}^r\sum_{i=1}^p (u_{i+1}\ldots u_pv_1\ldots
v_q\otimes w_{k+1}\ldots w_r)\omega_{n-2}\hfill\cr \hfill{\rm d}(u_1\ldots
u_{i-1}\otimes bw_1\ldots w_{k-1}) [{\rm d}(u_i\otimes 1)
{\rm d}(1\otimes w_k)-{\rm d}(1\otimes w_k){\rm d}(u_i\otimes
1)]\cr +(-1)^n\sum_{k=2}^r
(v_1\ldots v_q\otimes w_{k+1}\ldots w_r)\omega_{n-2} {\rm
d}(u_1\ldots u_p\otimes b){\rm d}(1\otimes w_1\ldots w_{k-1}){\rm
d}(1\otimes w_k)\hfill\cr +(-1)^n\sum_{j=1}^q (v_{j+1}\ldots
v_q\otimes 1)\omega_{n-2} {\rm d}(u_1\ldots u_p\otimes b){\rm
d}(v_1\ldots v_{j-1}\otimes w_1\ldots w_r){\rm d}(v_j\otimes
1)\hfill\cr}$$ (en particulier $h_0=0$ et dans $h_1$ seules les
deux derni\`eres lignes interviennent).
\end{lem}

\noindent{\sl Preuve.\/} Posons $C=A\otimes
B$, $C^e=C\otimes C^{op}$ et $N_n(C)=C^e\otimes\bar C^{\otimes
n}=C\otimes\bar C^{\otimes
n}\otimes C$. Alors $\Omega^n(C)=C\otimes_{C^e}N_n(C)$, or
$(N_n(C),b')$ est une r\'esolution libre du $C^e$-module $C$ (\cite{CE} 
p.~176). On construit une autre r\'esolution libre $M_*(C)$ de $C$,
{\sl de longueur 2,\/} en tensorisant deux r\'esolutions libres de
longueur 1 de $A$ et $B$. Celle pour $A$, par exemple, est donn\'ee
par $M_0(A)=A\otimes A=N_0(A)$, $M_1(A)=A\otimes V\otimes
A\subset N_1(A)$, avec une homotopie contractante $s$,
d\'efinie par~: $s$ est $A$-lin\'eaire \`a droite,
$s_{-1}(1)=[]$ et $s_0(v_1\ldots v_q[])=\sum_{j=1}^qv_1\ldots
v_{j-1}[v_j]v_{j+1}\ldots v_q$. La
r\'esolution $M_*(C)$ obtenue est donc d\'efinie par
$M_0(C)=C\otimes C$, $M_1(C)=C\otimes(W\oplus V)\otimes C$,
$M_2(C)=C\otimes(V\otimes W)\otimes C$, avec comme diff\'erentielle
$b'$, et comme homotopie contractante $\sigma $ d\'efinie par~:
$\sigma $ est $C$-lin\'eaire \`a droite, $\sigma_{-1}(1_C)=[]$,
$$\displaylines{\sigma _0((v_1\ldots v_q\otimes w_1\ldots w_r)[])=
\sum_{k=1}^r(1_A\otimes w_1\ldots w_{k-1})[w_k](v_1\ldots
v_q\otimes w_{k+1}\ldots w_r)\hfill\cr\hfill +
\sum_{j=1}^q(v_1\ldots v_{j-1}\otimes w_1\ldots
w_r)[v_j](v_{j+1}\ldots v_q\otimes 1_B),\cr
\sigma _1((v_1\ldots v_q\otimes b)[w])=\sum_{j=1}^q(v_1\ldots
v_{j-1}\otimes b)[v_j,w](v_{j+1}\ldots v_q\otimes
1_B)\cr}$$ et $\sigma _1(c[v])=0$.

Par la m\'ethode de \cite{CE} p.~76-77, on d\'efinit explicitement deux
morphismes ($C^e$-lin\'eaires) $F:M_*(C)\to N_*(C)$, $G:N_*(C)\to
M_*(C)$, au-dessus de ${\rm id}_C$, et une homotopie
$H_*:N_*(C)\to N_{*+1}(C)$ entre leur compos\'e $K:=FG$ et
l'identit\'e (de $N_*(C)$). Pr\'ecisons les r\'esultats
interm\'ediaires pour $F$ et $G$~:
$$\displaylines{F_0=G_0={\rm id}_{C\otimes C},\hfill
F_1([w])=[1_A\otimes w],\hfill
F_1([v])=[v\otimes 1_B],\cr
F_2([v,w])=[v\otimes
1_B,1_A\otimes w]-[1_A\otimes w,v\otimes 1_B],\cr
G_1([v_1\ldots v_q\otimes w_1\ldots w_r])=
\sum_{k=1}^r(1_A\otimes w_1\ldots w_{k-1})[w_k](v_1\ldots
v_q\otimes w_{k+1}\ldots w_r)+\hfill\cr\hfill
\sum_{j=1}^q(v_1\ldots v_{j-1}\otimes w_1\ldots
w_r)[v_j](v_{j+1}\ldots v_q\otimes 1_B),\cr
G_2([v_1\ldots v_q\otimes b,a\otimes
w_1\ldots w_r])=\hfill\cr\hfill \sum_{k=1}^r\sum_{j=1}^q(v_1\ldots
v_{j-1}\otimes bw_1\ldots w_{k-1})[v_j,w_k](v_{j+1}\ldots v_qa\otimes
w_{k+1}\ldots w_r).\cr}$$
Il suffit de tensoriser $H$ et $K$ par $1_C$
au-dessus de $C^e$ pour en d\'eduire $h$ et $k$.\cqfd

On calcule ensuite $\varphi $ \`a partir de $h$ comme expliqu\'e
juste avant le lemme ci-dessus, et l'on trouve~:

\begin{lem}\label{lem:barshbarAW}
L'application
$\varphi_n:\Omega^n(A\otimes B)\to\Omega^{n+1}(A\otimes B)$
suivante (naturelle en les alg\`ebres $A$ et $B$ et $\Lambda
$-repr\'esentable) v\'erifie~: $\bbar{sh}\
\bbar{AW}=1+b\varphi+\varphi b$.
$$\displaylines{\varphi _n:=\sum_{0<p\leq
r\leq n}(-1)^{n+r}\varphi _n^{r,p},\qquad{\rm avec}\hfill\cr
\varphi _n^{r,p}(\omega {\rm d}(a_1\otimes b_1)\ldots{\rm
d}(a_r\otimes b_r)):=\hfill\cr\hfill
(a_{p+1}\ldots a_r\otimes 1_B)\omega {\rm d}(b_1\ldots
b_p)Sh_{p,r-p}({\rm d}a_1,\ldots,{\rm d}a_p,{\rm
d}b_{p+1},\ldots,{\rm d}b_r),\cr}$$ o\`u $Sh_{p,q}$ d\'esigne la
somme de tous les $(p,q)$-shuffles multipli\'es par leur signature.
De plus, cette homotopie est sp\'eciale.
\end{lem}

\noindent{\sl Preuve.\/} Par r\'ecurrence, on calcule
$Z_n:=(\bbar{sh}\ \bbar{AW}-1-\varphi
_{n-1}b)(a_0\otimes b_0{\rm d}(a_1\otimes b_1)\ldots {\rm
d}(a_n\otimes b_n))$, $\Phi _n:=h_n(Z_n)$, et $\varphi _n$
(d\'eduit de $\Phi _n$ par universalit\'e), et l'on v\'erifie que
$\varphi _n$ est normalisable, que $b(Z_n)=0$, et que $k_n(Z_n)=0$ (m\^eme si
$n\leq 2$, valeurs pour lesquelles $k_n\neq 0$), donc $b(\Phi
_n)=(bh+hb)(Z_n)=(1-k)(Z_n)=Z_n$, d'o\`u $b\varphi _n=\bbar{sh}\ \bbar{AW}-1-\varphi
_{n-1}b$. On v\'erifie facilement que le $\varphi $ obtenu est
sp\'ecial, et $\Lambda $-repr\'esentable (corollaire \ref{coro:III.6}).\cqfd

Comparons notre $\bbar{sh}_\infty$ avec la coextension
\og~classique~\fg\ de $\bbar{sh}$ (\S \ref{subsection:shuffle}), c'est-\`a-dire notre
$\bbar{sh}_1=\varphi B\bbar{sh}$ avec le shuffle cyclique
$\bbar{sh}'$. Rappelons qu'on a introduit au \S \ref{subsection:formulation} deux op\'erateurs $\sigma
$, permettant d'exprimer la commutativit\'e de $sh$ par~: $\sigma sh\sigma
=sh$. Ces $\sigma $ commutent non seulement \`a $b$ mais aussi \`a $B$, et
sont norma\-lisables (leurs normalis\'es seront encore not\'es $\sigma $).

\begin{prop}
On a $\varphi
B\bbar{sh}=\sigma \bbar{sh}'\sigma $.
\end{prop}

\noindent{\sl Preuve.\/} Calculons
$\varphi_{p+q+1}B_{p+q}\bbar{sh}_{p,q}(a_0{\rm d}a_1\ldots{\rm
d}a_p\otimes b_0{\rm d}b_1\ldots{\rm d}b_q)$, c'est-\`a-dire
appliquons $(-1)^{p+q+1+r}\varphi _{p+q+1}^{r,p'}$, pour $0<p'\leq
r\leq p+q+1$, aux expressions de la forme $\varepsilon (\sigma
)(-1)^{(p+q)i}{\rm d}x_{p+q+1-i}\ldots{\rm d}x_{p+q}{\rm
d}x_0\ldots{\rm d}x_{p+q-i}$, pour $0\leq i\leq p+q$,
$x_0=a_0\otimes b_0$, et $(x_1,\ldots,x_{p+q})=$ un $(p,q)$-shuffle
$\sigma $ de $(a_1,\ldots, a_p,b_1,\ldots,b_q)$. Posons $n=p+q$,
$q'=r-p'$, et $$u_1=x_{n+1-i},\ldots, u_i=x_n, u_{i+1}=x_0,
u_{i+2}=x_1,\ldots, u_{n+1}=x_{n-i}.$$ Pour que $\varphi
_{n+1}^{r,p'}({\rm d}u_1\ldots{\rm d}u_{n+1})\neq 0$ il faut, vue
la forme de $\varphi _{n+1}^{r,p'}$ et puisque
$u_{i+1}=x_0=a_0\otimes b_0$ est le seul terme \og~mixte~\fg,
que~: 
\begin{itemize}
\item[$\star$] $(n+1)-r+1\leq i+1\leq (n+1)-q'$,
\item[$\star$] $u_{(n+1)-r+1},\ldots, u_i$ et
$u_{i+2},\ldots,u_{(n+1)-q'}$ soient des $a$ et
\item[$\star$] $u_{(n+1)-q'+1},\ldots,u_{n+1}$ soient des $b$.
\end{itemize}

De plus, comme $u_1,\ldots, u_{n+1}$ est le $i$-permut\'e
cyclique du shuffle $\sigma $, les $a$ et $b$ \'evoqu\'es sont
n\'ecessairement~:\\
$(u_{i+2},\ldots,u_{(n+1)-q'})=(a_1,\ldots, a_{n-q'-i})$ et
$(u_{n+2-r},\ldots, u_i)=(a_{2n+2-r-i},\ldots, a_p)$ (ce qui
implique $n-q'-i<n+p+2-r-i$), et\\
$(u_{n+2-q'},\ldots,u_{n+1})=(b_1,\ldots,b_{q'})$ (ce qui
implique $q'\leq q$).\\
Autrement dit, en posant $j=n-q'-i$ et $k=n+p+1-r-i$, le shuffle
est n\'ecessairement de la forme
$$\sigma =(a_1,\ldots,a_j,b_1,\ldots,
b_{q'},\sigma'(a_{j+1},\ldots,a_k,b_{q'+1},\ldots, b_q),
a_{k+1},\ldots, a_p),$$
o\`u $\sigma'$ est un $(k-j, q-q')$-shuffle.
Les signatures de $\sigma $ et $\sigma '$ sont donc reli\'ees par~:
$\varepsilon (\sigma )=(-1)^{q(p-j)+(q-q')(k-j)}\varepsilon
(\sigma ')$.

$\varphi _{n+1}^{r,p'}({\rm d}u_1\ldots{\rm d}u_{n+1})$ est alors
\'egal \`a~:
$$\displaylines{\sigma '({\rm d}a_{j+1},\ldots,{\rm d}a_k,{\rm
d}b_{q'+1},\ldots,{\rm d}b_q){\rm d}b_0Sh_{p',q'}({\rm
d}a_{k+1}\ldots{\rm d}a_p{\rm d}a_0\ldots{\rm d}a_j,{\rm
d}b_1\ldots{\rm d}b_q')\cr =\sigma ''({\rm d}b_0,\ldots,{\rm
d}b_q,{\rm d}a_0,\ldots,{\rm d}a_p),\cr}$$ o\`u $\sigma ''$ est un
$(q+1,p+1)$-shuffle cyclique (cf. \S \ref{subsection:shuffle}).
Les signatures de $\sigma' $ et $\sigma ''$ sont donc reli\'ees par~:
$\varepsilon (\sigma ''
)=(-1)^{q(q-q')+p(p-j)+q'(p+1)+(k-j)(q-q'+1)}\varepsilon (\sigma
')$.

Les conditions sur $(i,r,q')$~:
$0\leq q'(\leq r-1\leq n)$, $(0\leq i\leq n)$, $n+1-r\leq i\leq
n-q'$, $n-q'-i<n+p+2-r-i$, $q'\leq q$
sont redondantes et \'equivalent simplement \`a~:
$0\leq q'\leq q$, $q'+i\leq n$, $r\leq q'+p+2$, $n< r+i$.
Elles \'equivalent donc \`a~:
$0\leq q'\leq q, 0\leq j\leq k\leq p$.
Les $\sigma ''$ obtenus sont donc {\sl tous\/} les
$(q+1,p+1)$-shuffles cycliques possibles (une fois chacun) de ${\rm
d}b_0,\ldots,{\rm d}b_q,{\rm d}a_0,\ldots,{\rm d}a_p$.

Compte tenu des liens entre les signatures de $\sigma ,\sigma
',\sigma ''$, le coefficient de $\sigma ''$ dans $\varphi Bsh(a_0{\rm
d}a_1\ldots{\rm d}a_p\otimes b_0{\rm d}b_1\ldots{\rm d}b_q)$ est
$(-1)^{ni+n+1+r}\varepsilon (\sigma )=(-1)^{qp+q}\varepsilon
(\sigma '')$.
On aboutit donc \`a~: $$\displaylines{\varphi Bsh(a_0{\rm
d}a_1\ldots{\rm d}a_p\otimes b_0{\rm d}b_1\ldots{\rm
d}b_q)=\hfill\cr\hfill (-1)^{qp+q} \sum_{\sigma '' \in S'_{q+1,p+1}}
\varepsilon (\sigma '' )\sigma ''({\rm d}b_0,\ldots,{\rm
d}b_q,{\rm d}a_0,\ldots,{\rm d}a_p),\cr}$$ ce qui (\S \ref{subsection:formulation} et \ref{subsection:shuffle}) correspond
bien \`a la d\'efinition de $\sigma \bbar{sh}'\sigma $.\cqfd

\noindent{\sl Remarques.\/}
\begin{itemize}
\item Contrairement \`a $\bbar{sh}$, $\bbar{AW}$ n'est pas commutatif. On
avait donc en fait le choix, dans l'explicitation du th\'eor\`eme \ref{thm:perturbcoext}, entre
deux quasi-inverses naturels (et homotopies associ\'ees) pour $\bbar {sh}$,
donnant deux coextensions diff\'erentes de $\bbar{sh}$, conjugu\'ees l'une de
l'autre par $\sigma $~: le quasi-inverse que nous avons choisi,
$(\bbar{AW},\varphi )$, a donn\'e une coextension par $\varphi B\bbar{sh}$,
l'autre, $(\sigma \bbar{AW}\sigma ,\sigma \varphi\sigma  )$, aurait donn\'e une
coextension par $\bbar{sh}'$.
\item Le th\'eor\`eme IV.3 et la proposition IV.6
fournissent donc une autre preuve de \cite{L} 4.3.3 (en tenant compte de la
rectification signal\'ee au \S \ref{subsection:shuffle}), c'est-\`a-dire de
$[B,\bbar{sh}]+[b,\bbar{sh}']=0$.
\end{itemize}

\begin{prop}\label{prop:hassoc}
$\bbar{sh}_\infty$ est
associatif \`a homotopie pr\`es, et $\bbar{AW}_\infty$ est
coassociatif \`a homotopie pr\`es.
\end{prop}

\noindent{\sl Preuve.\/} Il suffit \'evidemment de prouver
l'associativit\'e (\`a homotopie pr\`es) de $\bbar{sh}_\infty$.
Par associativit\'e de $\bbar{sh}$ et par des raisonnements sur
les degr\'es (comme dans la remarque \ref{rem:normalisable}), il suffit en fait de
trouver (pour tous $p,q,r\geq 0$) un
$\psi:\Omega^p(A)\otimes\Omega^q(B)\otimes\Omega^r(C)\to\Omega^{p+q+r+3}(A\otimes
B\otimes C)$ tel que $$\bbar{sh}(\bbar{sh}_1\otimes
1)+\bbar{sh}_1(\bbar{sh}\otimes
1)-\bbar{sh}(1\otimes\bbar{sh}_1)-\bbar{sh}_1(1\otimes\bbar{sh})=b\psi+\psi
b.$$ Le lecteur est invit\'e, \`a titre d'exercice, \`a utiliser, pour
construire un tel $\psi$, la m\^eme m\'ethode que celle pr\'esent\'ee pour
construire $\varphi $ (lemmes \ref{lem:homotopie} et \ref{lem:barshbarAW}). Mais la proposition pr\'ec\'edente
permet d'\'eviter ces calculs, en utilisant les op\'erateurs $B_2, B_3$ de
[GJ2]. En effet, le $B_2$ de Getzler-Jones n'est autre que $\bbar{sh}'$ au
signe pr\`es~: $\bbar{sh}'_{p,q}=(-1)^pB_2$, donc l'\'equation de
[GJ2], lemma 4.3, $$\displaylines{-bB_3(\alpha ,\beta ,\gamma
)=(-1)^{\varepsilon _0}B_3(b\alpha ,\beta ,\gamma
)+(-1)^{\varepsilon _1}B_3(\alpha ,b\beta ,\gamma
)+(-1)^{\varepsilon _2}B_3(\alpha ,\beta ,b\gamma
)\hfill\cr\hfill+\alpha\star B_2(\beta ,\gamma )+(-1)^{\varepsilon
_1}B_2(\alpha \star\beta ,\gamma )+(-1)^{\varepsilon _2}B_2(\alpha
,\beta \star \gamma )-(-1)^{\varepsilon _2}B_2(\alpha ,\beta
)\star\gamma\cr}$$ avec $\alpha ,\beta ,\gamma $ de degr\'es $p,q,r$
et $\varepsilon _0=0,\varepsilon _1=p-1,\varepsilon _2=p+q$,
devient, en posant $\psi(\alpha ,\beta ,\gamma )=(-1)^qB_3(\alpha
,\beta ,\gamma )$~: $$\bbar{sh}(\bbar{sh}'\otimes
1)+\bbar{sh}'(\bbar{sh}\otimes
1)-\bbar{sh}(1\otimes\bbar{sh}')-\bbar{sh}'(1\otimes\bbar{sh})=b\psi+\psi
b.\qquad\Box$$

On pourrait esp\'erer, puisque $\bbar{sh}$ est commutatif, que
$\bbar{sh}_\infty$ le soit \`a homotopie (naturelle) pr\`es,
c'est-\`a-dire que les deux coextensions $\bbar{sh}_\infty$ et
$\sigma \bbar{sh}_\infty\sigma $ de $\bbar{sh}$ soient
(naturellement) homotopes. C'est ce qu'affirme le th\'eor\`eme
2.3 de \cite{HJ}, mais nous allons prouver que {\sl ce
th\'eor\`eme est faux.\/} (Les coextensions de $\bbar{sh}$ \`a homotopie pr\`es
forment en r\'ealit\'e une droite affine, de direction $H_2(\Lambda
[0]\times\Lambda [0])=k[B\otimes B]$).

\begin{prop}
$\bbar{sh}_\infty$ et
$\sigma \bbar{sh}_\infty\sigma $ ne sont pas (naturellement)
homotopes.
\end{prop}

\noindent{\sl Preuve.\/} Montrons que
$\psi_\infty:=\bbar{AW}_\infty(\sigma \bbar{sh}_\infty\sigma-\bbar{sh}_\infty)$
n'est pas (naturellement) homotope \`a $0$. Si elle l'\'etait, il existerait des
applications (naturelles)\\
$h^{(i)}_k:(\Omega(A)\otimes\Omega(B))_k\to
(\Omega(A)\otimes\Omega(B))_{k+2i+1}$ (pour $i,k\geq 0$) telles que
$$\psi^{(i)}_k=bh^{(i)}_k+h^{(i)}_{k-1}b+Bh^{(i-1)}_k+h^{(i-1)}_{k+1}B.$$ Mais
pour des raisons de degr\'es (remarque \ref{rem:normalisable}), de tels $h^{(i)}_k$ seraient
nuls pour $i\geq 1$, donc il ne resterait que les $h_k:=h^{(0)}_k$. D'autre
part, $\sigma \bbar{sh}\sigma =\bbar{sh}$, donc
$\psi^{(0)}=0$ et $$\psi^{(1)}=\bbar{AW}^{(1)}.0+\bbar{AW}(\sigma
\bbar{sh}^{(1)}\sigma -\bbar{sh}^{(1)})=\bbar{AW}\ \bbar{sh}'-(\bbar{AW}\varphi)
B\bbar{sh}=\bbar{AW}\ \bbar{sh}',$$ en
particulier $\psi^{(1)}_0(a\otimes b)=\bbar{AW}\  \bbar{sh}'(a\otimes
b)=\bbar{AW}({\rm d}a{\rm d}b)={\rm d}a\otimes{\rm d}b$. On aurait donc~:
$$bh_k=-h_{k-1}b\qquad {\rm et}\qquad {\rm d}a\otimes{\rm
d}b=(Bh_0+h_1B)(a\otimes b).$$ Or $h_0:=h_{0,0}$ et $h_1 :=(h_{1,0},h_{0,1})$
sont a priori de la forme suivante~: $$\displaylines{h_{0,0}(a\otimes
b)=\alpha {\rm d}a\otimes b+\beta a\otimes{\rm d}b\hfill\cr h_{1,0}(a{\rm
d}a'\otimes b)= (\alpha _1{\rm d}a{\rm d}a'+\alpha _2{\rm d}a'{\rm d}a)\otimes
b+\hfill\cr\hfill (\beta _1 a{\rm d}a'+\beta _2a'{\rm d}a+\gamma _1{\rm
d}(aa')+\gamma _2{\rm d}(a'a))\otimes {\rm d}b\cr h_{0,1}(a\otimes b{\rm
d}b')= a\otimes(\alpha _3{\rm d}b{\rm d}b'+\alpha _4{\rm d}b'{\rm d}b)+\hfill\cr\hfill 
{\rm d}a\otimes(\beta _3 b{\rm d}b'+\beta _4b'{\rm d}b+\gamma _3{\rm
d}(bb')+\gamma _4{\rm d}(b'b)),\cr}$$ avec, pour que
$bh_{1,0}=-h_{0,0}(b_1\otimes 1_0)$ et
$bh_{0,1}=-h_{0,0}(1_0\otimes b_1)$~:
$$-\alpha _2=\alpha _1=\alpha,\qquad -\alpha_4=\alpha _3=\beta,\qquad\beta
_1-\beta _2=-\beta,\qquad \beta _3-\beta_4=+\alpha.$$ D'autre part, pour que
$-h_{1,0}(b_2\otimes 1_0)$ et $-h_{0,1}(1_0\otimes b_2)$ puissent se mettre
respectivement sous la forme $bh_{2,0}$ et $bh_{0,2}$, il faut (et il suffit)
que $$-\alpha _2=\alpha _1,\qquad -\alpha_4=\alpha _3,\qquad\gamma _1+\gamma
_2=-\beta _2,\qquad\gamma _3+\gamma _4=-\beta _4.$$ On aurait
alors~: $${\rm d}a\otimes{\rm d}b=(Bh_0+h_1B)(a\otimes b)= (-\alpha +\beta
+\beta _1+\gamma _1+\gamma _2+\beta _3+\gamma _3+\gamma _4){\rm d}a\otimes{\rm
d}b=0,$$ d'o\`u la contradiction.\cqfd

\begin{rem}\label{rem:calculinit}
Calcul des \og~initialisations~\fg\  
de $\bbar{sh}$ et $\bbar{AW}$ (cf. \S \ref{subsection:Hood-Jones} et \ref{subsection:Kassel}). Posons
$F^{(1)}=\sigma \varphi \sigma B$ (de mani\`ere \`a avoir
$\bbar{sh}'=F^{(1)}\bbar{sh}$). Alors $F^{(1)}_0(a\otimes b)={\rm
d}a{\rm d}b$, ce qui correspond bien \`a l'initialisation donn\'ee
dans \cite{K2} p.~208, qu'il faut compl\'eter par $$\displaylines{F^{(1)}_1(a\otimes
b{\rm d}(a'\otimes b'))={\rm d}(a'a){\rm
d}bdb'-{\rm d}(a'a){\rm d}b'{\rm d}b+{\rm d}b'{\rm d}(a'a){\rm
d}b\hfill\cr\hfill -{\rm d}a{\rm d}a'{\rm d}(bb')+{\rm d}a'{\rm d}a{\rm
d}(bb')+{\rm d}a{\rm d}(bb'){\rm d}a'.\cr}$$ A partir de $F^{(0)}$ et $F^{(1)}$
on (re)trouve~:
$$\displaylines{\bbar{sh}'_{0,0}(a\otimes b)={\rm d}a{\rm d}b,\cr
\bbar{sh}'_{1,0}(a{\rm d}a'\otimes b)=
-{\rm d}a{\rm d}a'{\rm d}b+{\rm d}a'{\rm d}a{\rm d}b+
{\rm d}a{\rm d}b{\rm d}a',\cr
\bbar{sh}'_{0,1}(a\otimes b{\rm d}b')=
{\rm d}a{\rm d}b{\rm d}b'-{\rm d}a{\rm d}b'{\rm d}b+
{\rm d}b'{\rm d}a{\rm d}b,\cr}$$
ce qui est conforme \`a notre d\'efinition du shuffle cyclique (\S
\ref{subsection:shuffle}) et correspond bien \`a l'initialisation donn\'ee (sous forme de
$\Lambda$-repr\'esentants) dans \cite{HJ} p.~371. Par contre la formule
correspondante de \cite{HJ} p.~373 pour $\bbar{AW}^{(1)}$ est fausse. On peut la
remplacer par notre $\bbar{AW}^{(1)}:=\bbar{AW}B\varphi$, qui
donne
$$\displaylines{\bbar{AW}^{(1)}_0(a\otimes b)=0,\cr
\bbar{AW}^{(1)}_1(a\otimes b{\rm d}(a'\otimes b'))=
a{\rm d}a'\otimes {\rm d}b{\rm d}b'+
{\rm d}a'{\rm d}a\otimes b{\rm d}b'\cr}$$
(ou, en termes de $\Lambda $-repr\'esentants~:
$\bbar{AW}^{(1)}_1=1\otimes t_2s_1+t_2^2s_0\otimes 1$).
\end{rem}

\begin{rem}\label{rem:DeRham}
Soient $A$ une alg\`ebre commutative et
$(\Lambda (A),{\rm d})$ son complexe de De Rham commutatif. L'application
$\pi :\omega\mapsto{\bar\omega \over n!}$ est un morphisme de
$(\Omega(A),B,b)$ sur $(\Lambda (A),{\rm d}, 0)$ (\cite{CQ2}, \S 13). Pour $A$ et $B$
commutatives, le produit $\bbar{sh}_\infty$ passe aux quotients, et est
compatible avec le produit usuel. Plus pr\'ecis\'ement~: $$\pi\bbar{sh}_\infty
(\omega _p\otimes\omega _q)=\pi (\omega _p)\wedge\pi (\omega
_q)\pm{(-1)^p\over 2}{\rm d}\pi (\omega _p)\wedge{\rm d}\pi (\omega _q),$$
avec $\pm=+$ si l'on a choisi $\bbar{sh}'$ et $\pm=-$ si l'on a choisi
$\sigma \bbar{sh}'\sigma $ (pour co\'etendre $\bbar{sh}$), mais peu importe,
puisqu'on v\'erifie facilement que $x_p\otimes y_q\mapsto{(-1)^p\over 2}{\rm
d}x_p\wedge{\rm d}y_q$ est homotope \`a $0$.
\end{rem}

Il resterait, pour compl\'eter le tableau, \`a expliciter le lien
avec le produit en cohomologie cyclique d\'efini par Connes.

\section{D\'enormalisation d'une coextension.}\label{section:denormalisation}

Les deux \'equations $\bbar{AW}\ \bbar{sh}=1$,
$\bbar{AW}B\bbar{sh}=B$ qui ont permis, au paragraphe
pr\'ec\'edent, de construire des coextensions $\bbar{AW}_\infty$,
$\bbar{sh}_\infty$ de $\bbar{AW}$ et $\bbar{sh}$, ne sont
satisfaites que dans les complexes normalis\'es. Donc la m\'ethode
du paragraphe pr\'ec\'edent n'est pas applicable \`a $AW$ et $sh$.
Par contre, nous allons \og~relever\fg\
$\bbar{AW}_\infty$, $\bbar{sh}_\infty$ en des coextensions
$AW_\infty$, $sh_\infty$ des $b$-morphismes $AW$, $sh$ (corollaire
V.2). Plus g\'en\'eralement, la proposition suivante montre comment
relever une coextension $g_\infty$ d'un $b$-morphisme
$g_0$ \`a valeurs dans un complexe normalis\'e, {\sl sans
supposer a priori que $g_0$ se rel\`eve en un $b$-morphisme\/} (\`a valeurs
dans le complexe non normalis\'e). On peut remarquer (en vue d'un \'enonc\'e
plus g\'en\'eral) que la seule propri\'et\'e de la normalisation $j:Q_*\to \bar
Q_*$ qu'on utilise est que c'est une r\'etraction par d\'eformation {\sl qui
commute \`a $B$.\/} Dans la proposition et le corollaire suivants, $\varphi $
d\'esignera une homotopie \og~sp\'eciale\fg\ associ\'ee \`a $j$ (cf. \S \ref{subsection:perturbnorma}). 

\begin{prop}
Soient $f^{(k)}:P_*\to Q_{*+2k}$
des applications telles que les $g^{(k)}:=jf^{(k)}:P_*\to \bar
Q_{*+2k}$ v\'erifient $[b,g^{(k)}]+[B,g^{(k-1)}]=0$.
Alors les $h^{(k)}:P_*\to Q_{*+2k}$ d\'efinies par
$$h^{(k)}_n:=f^{(k)}_n+\varphi
(bf^{(k)}_n-h^{(k)}_{n-1}b+[B,h^{(k-1)}]_n)$$
v\'erifient~: $jh^{(k)}=g^{(k)}$ et $[b, h^{(k)}]+[B, h^{(k-1)}]=0$.
\end{prop}

\noindent{\sl Preuve.\/} On pose bien s\^ur $h^{(k)}_n=0$
si $k$ ou $n$ $<0$, ce qui initialise cette d\'efinition par r\'ecurrence (en
particulier, on aura $h_0^{(0)}=f_0^{(0)}$). Puisque $j\varphi =0$,
on a d\'ej\`a $jh=jf=g$. Posons alors $x^{(k)}_n:=([b,
h^{(k)}]+[B,h^{(k-1)}])_n$ (dont on veut montrer qu'il est nul) et
$y:=x-bh+bf$. On a donc~: $(1+b\varphi +\varphi b)y=0$ et $by=bx$, d'o\`u
$bh=b(f+\varphi y)=bf-(1+\varphi b) y=bf-y-\varphi bx=bh-x-\varphi bx$,
donc $(1+\varphi b)x=0$. Or
$$\displaylines{bx^{(k)}=-bh^{(k)}b-Bbh^{(k-1)}-bh^{(k-1)}B=\cr
([B,h^{(k-1)}]-x^{(k)})b
-B(x^{(k-1)}+h^{(k-1)}b+h^{(k-2)}B)
-(x^{(k-1)}+h^{(k-1)}b-Bh^{(k-2)})B\cr
=-x^{(k)}b-Bx^{(k-1)}-x^{(k-1)}B.\cr}$$
On a donc $x^{(k)}_n=-\varphi
bx^{(k)}_n=\varphi (x^{(k)}_{n-1}b+Bx^{(k-1)}_n+x^{(k-1)}_{n+1}B)$, d'o\`u
l'on d\'eduit (par r\'ecurrence sur $(k,n)$) que $x^{(k)}_n=0$.\cqfd

\begin{coro}
Soient $f^{(k)}:P_*\to Q_{*+2k}$
des applications telles que les $g^{(k)}:=jf^{(k)}:P_*\to \bar
Q_{*+2k}$ v\'erifient $[b,g^{(k)}]+[B,g^{(k-1)}]=0$, et telles que
$[b,f^{(0)}]=0$. Alors les $h^{(k)}:P_*\to Q_{*+2k}$ d\'efinies par
$$h^{(k)}:=f^{(k)}+\varphi \sum_{i=1}^k(B\varphi
)^{k-i}([b,f^{(i)}]+[B,f^{(i-1)}])$$ v\'erifient~: $jh^{(k)}=g^{(k)}$, $[b,
h^{(k)}]+[B, h^{(k-1)}]=0$, et $h^{(0)}=f^{(0)}$.
\end{coro}

\noindent{\sl Preuve.\/} Il suffit de montrer que gr\^ace \`a l'hypoth\`ese
suppl\'ementaire sur $f^{(0)}$, ces $h$ co\"\i ncident avec les
pr\'ec\'edents. La v\'erification (par r\'ecurrence sur $(k,n)$) est
imm\'ediate (en utilisant $\varphi ^2=0$).\cqfd

\begin{rem}\label{rem:V.3}
Si $g_\infty$ est de longueur 2, c'est-\`a-dire
si $f^{(k)}$ peut \^etre suppos\'ee nulle pour $k\geq 2$ (ce qui sera le cas
pour $g_\infty=\bbar{sh}_\infty$), son \og~d\'enormalis\'e\fg\ $h_\infty$ ne le
restera pas~: on aura seulement, pour $k>0$,
$h^{(k)}=f^{(k)}+\varphi (B\varphi )^{k-1}([b,f^{(1)}]+[B,f^{(0)}])$, d'o\`u
l'int\'er\^et du \S \ref{section:coext-long-2}.
\end{rem}

\begin{rem}
Le corollaire ci-dessus peut en fait se
d\'eduire directement du lemme de perturbation (\S \ref{subsection:perturbnorma}) qui, appliqu\'e
\`a $(j,\varphi )$, donne $b_\infty=b+B$ et $j_\infty=j$. En effet, la
d\'efinition des $h^{(k)}$ s'\'ecrit $h_\infty=f_\infty+\varphi
_\infty[b_\infty,f_\infty]$ et l'hypoth\`ese sur les
$f^{(k)}$ donne $(1+b_\infty\varphi _\infty+\varphi _\infty
b_\infty)[b_\infty,f_\infty]=i_\infty j[b_\infty,f_\infty]=0$ donc
$[b_\infty ,h_\infty ]=[b_\infty ,f_\infty ]+b_\infty\varphi_\infty [b_\infty,
f_\infty]-\varphi_\infty[b_\infty,f_\infty]b_\infty=-\varphi_\infty b_\infty
[b_\infty,f_\infty]-\varphi_\infty [b_\infty,f_\infty]b_\infty=0$.
De plus, d'apr\`es notre reformulation (lemme \ref{lem:perturbation}), il est en
r\'ealit\'e inutile ici de supposer que $\varphi $ est sp\'eciale. On a juste
besoin (pour que $jh_\infty=jf_\infty$) que $j\varphi =0$, mais pas que
$\varphi ^2$ ni $\varphi i$ soient nulles.
\end{rem}

\begin{rem}\label{rem:V.5}
On peut appliquer le corollaire
ci-dessus pour \og~d\'enormaliser\fg\ $\bbar{AW}$ et $\bbar{sh}$. Il n'est
pas n\'ecessaire pour cela d'expliciter $i$ et $j$~: il suffit de $\varphi $.
D'apr\`es la fin de la remarque pr\'ec\'edente, on peut choisir (en
s'inspirant de \cite{M} p.~94--95)
$$\displaylines{\varphi _n=-s_0+s_1(1-s_0d_0)-s_2(1-s_1d_1)(1-s_0d_0)+\ldots
\hfill\cr
\hfill+(-1)^{n+1}s_n(1-s_{n-1}d_{n-1})\ldots (1-s_0d_0).\cr}$$
On en d\'eduit en particulier une d\'enormalisation des initialisations de la
remarque \ref{rem:calculinit}~:
$$\displaylines{sh^{(1)}_{0,0}(a\otimes b)=(1,a,b)+(ab,1,1),\hfill\cr
sh^{(1)}_{1,0}((a,a')\otimes b)=
-(1,a,a',b)+(1,a',a,b)+(1,a,b,a')+\hfill\cr\hfill
(a',1,1,ab)-(a'b,1,1,a)+(ab,1,1,a')-(ab,1,a',1)
+(ab,a',1,1)-(aa'b,1,1,1),\cr
sh^{(1)}_{0,1}(a\otimes
(b,b'))=(1,a,b,b')-(1,a,b',b)+(1,b',a,b)+\hfill\cr\hfill
(b',1,1,ab)+(ab,1,1,b')-(ab',1,1,b)-(ab,1,b',1)+(ab,b',1,1)-(abb',1,1,1).\cr}$$
($AW_0^{(1)}=0$, et on laisse au lecteur le soin de calculer $AW_1^{(1)}$).
\end{rem}

\section{Produit et coproduit en homologie cyclique
enti\`ere.}\label{section:entiere}

Une $b+B$-cocha\^\i ne $(\varphi _{2n})_{n\in{\Bbb N}}$ (resp.
$(\varphi _{2n+1})_{n\in{\Bbb N}}$) sur une alg\`ebre de Banach $A$
est dite enti\`ere si et seulement si $\forall \mu  >0$,
$\sum_n \|\varphi_{2n}\| \mu ^n\ n!<\infty$
(resp. $\sum_n\|\varphi_{2n+1}\| \mu ^n\ n!<\infty$),
ou encore si et seulement si $\forall \mu >0$, la suite $\|\varphi
_{2n}\| \mu ^n\ n!$ (resp. $\|\varphi _{2n+1}\| \mu ^n\ n!$) est
born\'ee . (La d\'efinition originelle de la cohomologie cyclique
enti\`ere -- \cite{C1} -- \'etait en termes de $d_1+d_2$-cocha\^\i nes,
avec $d_1:=(n+1)b:C^n\to C^{n+1}$ et $d_2:={1\over n}B:C^n\to C^{n-1}$, mais
un tableau de conversion est fourni dans \cite{C2}, p.~371). Donc en unifiant les
cas pair et impair~: $(\varphi _n)_{n\in{\Bbb N}}$ est enti\`ere si et
seulement si $\forall r>0$, $\|\varphi _n\|  r^{-n}\sqrt{n!}$ est born\'ee.
{\sl En homologie,\/} le $b+B$ complexe entier (normalis\'e), est donc (\cite{GS})
$\Omega_\varepsilon (A):=$ le {\sl syst\`eme\/} inductif, pour $r\to 0^+$, des
$\Omega_r(A)$, d\'efinis comme compl\'etions de $\Omega(A)$ pour les normes
$$\|\sum_{n=0}^\infty \omega_n\|_r=\sum_{n=0}^\infty
{r^n\|\omega_n\|\over\sqrt{n!}}$$ (On v\'erifie sans peine que $b$ et $B$ sont
bien des endomorphismes du Ind-objet $\Omega_\varepsilon (A)$). J'insiste sur
le choix d'arr\^eter la construction aux Ind-objets au lieu de passer \`a la
limite inductive, qui n'est pas compatible au produit tensoriel projectif.

La preuve du th\'eor\`eme ci-dessous
est essentiellement combinatoire, donc\break s'adapte verbatim aux
variantes naturelles de cette d\'efinition pour des alg\`ebres
localement convexes (\cite{C2} p.~370), et aux th\'eories cycliques
(p\'eriodique, enti\`ere, asymptotique) \'etudi\'ees par Puschnigg
\cite{P1}.

\begin{thm}
Le produit et le coproduit
$\bbar{sh}_\infty,\bbar{AW}_\infty$ d\'efinissent des
quasi-isomorphismes, quasi-inverses mutuels,
$$\displaylines{\bbar{sh}_\varepsilon:\Omega_\varepsilon (A\otimes
B)\to\Omega_\varepsilon (A)\otimes\Omega_\varepsilon (B),\cr
\bbar{AW}_\varepsilon :\Omega_\varepsilon
(A)\otimes\Omega_\varepsilon (B)\to\Omega_\varepsilon (A\otimes
B).\cr}$$
\end{thm}

Il est facile de v\'erifier que $\varphi$, $sh$, $AW$
sont continues. Par exemple, la continuit\'e de $\varphi $ vient du
fait que $\varphi _n$ est une somme de $\sum_{0<p\leq r\leq
n}C^p_r=2(2^n-1)-n$ termes. Vues les formules du lemme \ref{lem:perturbation}, le seul probl\`eme
est donc de prouver la continuit\'e de l'application $\sum_{k>0}(\varphi
B)^k$, de $\Omega_\varepsilon (A\otimes B)$ dans lui-m\^eme. Dans cette
s\'erie, appliqu\'ee \`a un \'el\'ement $\omega_n$ de degr\'e $n$, les
termes de degr\'es $>2n+3$ sont nuls (remarque \ref{rem:normalisable}) donc la somme est en
r\'ealit\'e finie (pour $k$ de $1$ \`a $\lfloor{n+3\over
2}\rfloor$), mais la majoration pr\'ec\'edente du nombre de termes dans
$\varphi _n$ est trop grossi\`ere, car elle ne donne, pour le nombre de termes
de $(\varphi B)^k(\omega _n)$, qu'un majorant de l'ordre de $2^{k(n+k)}$, en
particulier, pour $k=\lfloor{n+3\over
2}\rfloor$, de l'ordre de $2^{(3n^2/4)}$, qui n'est pas major\'e par
une expression de la forme $C^n\sqrt{(2n+3)!\over n!}$.

Pour affiner la majoration de $\varphi B$, introduisons quelques
remarques et notations. Rappelons (lemme \ref{lem:barshbarAW}) que sur
$\Omega^n(A\otimes B)$, $\varphi _n:=\sum_{0<p\leq r\leq n}(-1)^{n+r}\varphi
_n^{r,p}$ avec, si $\omega \in \Omega^{n-r}(A\otimes B)$ et $m_i=a_i\otimes
b_i$,  $$\displaylines{\varphi _n^{r,p}(\omega {\rm d}m_1\ldots{\rm
d}m_r):=\hfill\cr\hfill (a_{p+1}\ldots a_r\otimes 1_B)\omega {\rm d}(b_1\ldots
b_p)Sh_{p,r-p}({\rm d}a_1,\ldots,{\rm d}a_p,{\rm d}b_{p+1},\ldots,{\rm
d}b_r).\cr}$$ Pour que cette expression soit non nulle il faut que
\begin{itemize}
\item[$\star$] $a_1,\ldots, a_p$ soient diff\'erents de $1$ et l'un au moins des $m_1,\ldots,
m_p$ soit \og~mixte\fg~, c'est-\`a-dire v\'erifie \'egalement $b_i\neq 1$, et
\item[$\star$] $b_{p+1},\ldots, b_r$ soient diff\'erents de $1$.
\end{itemize}

On en d\'eduit la propri\'et\'e suivante. Appelons \og~mot autoris\'e \fg~
un\break $m_0{\rm d}m_1\ldots{\rm d}m_n$ o\`u $(m_0,\ldots,m_n)$ est
constitu\'e (\`a permutation circulaire pr\`es) d'une succession de \og~blocs
vivants\fg\ (c'est-\`a-dire de suites $\alpha ^p\mu ^q\beta ^r$ de $p+q+r$
termes, les $p$ premiers de la forme $a\otimes 1$, les $q$ suivants ($q>0$)
mixtes, et les $r$ derniers de la forme $1\otimes b$), alternant avec des
\og~blocs inertes\fg\ (c'est \`a dire vides, ou compos\'es uniquements de
termes $\alpha $ et $\beta $, avec un $\alpha $ au d\'ebut et un $\beta $ \`a
la fin). Alors $\varphi B$ transforme un mot autoris\'e en une somme
(sign\'ee) de mots autoris\'es, obtenue en disloquant (un par un) chaque bloc
vivant et en faisant la somme.

Pr\'ecisons la dislocation d'un tel bloc $\alpha^p\mu ^q\beta ^r$~:
\begin{itemize}
\item on choisit dans $\alpha^p\mu ^q\beta ^r$ un sous-bloc \`a disloquer,
c'est-\`a-dire le nombre $(i)$ de $\alpha $ et $(k)$ de $\mu $ \'ecart\'es \`a
gauche, et le nombre $(j)$ de $\beta $ et $(m-k)$ de $\mu $ \'ecart\'es \`a
droite, avec $0\leq i\leq p$, $0\leq j\leq r$,$0\leq k\leq m<q$, et
$k>0\Rightarrow i=p$, $k<m\Rightarrow j=r$
\item puis on am\`ene ce sous-bloc $\alpha ^{p-i}\mu ^{q-m}\beta ^{r-j}$
en queue (par permutation circulaire due \`a
$B$) et on fait agir les $\varphi ^{r',p'}$ de
mani\`ere \`a remplacer (\`a permutation circulaire pr\`es dans le mot, et
aux signes pr\`es) ce sous-bloc par~: $\beta Sh_{p',r'-p'}(\alpha ^{p'}\beta
^{r'-p'})\alpha $
\item on \og~gomme\fg\ , dans chaque mot de la somme obtenue, des $\alpha
,\beta $ devenus \og~inertes\fg~.
\end{itemize}

Associons, \`a tout mot (sign\'e) sur l'alphabet $\alpha
,\beta ,\mu $, le mon\^ome $X^qY^s$, o\`u $q$ d\'esigne le nombre des $\mu
$, et $s$ la longueur du mot. Par exemple, au bloc vivant $\alpha^p\mu ^q\beta
^r$ \`a disloquer par $\varphi B$ est associ\'e le mon\^ome $X^qY^{p+q+r}$, et
au r\'esultat de cette dislocation (apr\`es gommage de certains $\alpha ,\beta
$ inertes) est associ\'e un polyn\^ome \`a coefficients positifs.
L'int\'er\^et de ce polyn\^ome est que sa valeur pour $X=Y=1$ donne le nombre
de termes de la somme obtenue par dislocation. Nous allons
majorer ce polyn\^ome (au sens~: majoration coefficient par coefficient), par
un polyn\^ome explicite ne d\'ependant que de $(q,p+q+r)$, que nous noterons
donc $(B\varphi )(X^qY^{p+q+r})$.

\begin{lem}
Il existe des constantes
universelles $K,L,M$  telles que~: $$(\varphi B)(X^q Y^s)\leq Kq\sum_{0\leq
m\leq n, m<q,n\leq s+2}
 L^{q-m}M^{s-n}X^mY^n.$$
\end{lem}

\noindent{\sl Preuve.\/} Lors de la dislocation du sous-bloc
$\alpha ^{p-i}\mu ^{q-m}\beta ^{r-j}$, de longueur $r':=p+q+r-i-j-m$,
le bloc $\alpha ^p\mu ^q\beta ^r$ est remplac\'e par une somme (sign\'ee) de
blocs, de la forme $\sum_{p'}\alpha^i\mu ^k\beta Sh_{p',r'-p'}(\alpha
^{p'}\beta ^{r'-p'})\alpha\mu ^{m-k}\beta ^j$. Dans une telle expression,
les seules lettres \'eventuellement vivantes d'un $(p',q')$-shuffle sont ses
$\beta $ les plus \`a gauche et ses $\alpha $ les plus \`a droite, et l'on
peut m\^eme effectuer des \og~gommages\fg\ suppl\'ementaires, lorsque $i>0$ et
$k=0$, ou lorsque $j>0$ et $k=m$. Il convient donc de distinguer les 3 cas
$$(k=i=0)\ {\rm ou}\ (0<i<p\ {\rm et}\ k=0)\ {\rm ou}\ (i=p\ {\rm et}\
k>0),$$ et dans chaque cas, les 3 sous-cas
$$(m-k=j=0)\ {\rm ou}\ (0<j<r\ {\rm et}\ k=m)\ {\rm ou}\ (j=r\ {\rm et}\
k<m),$$
soit au total 9 cas.
Notons $sh(\beta ,\alpha )$ ce qui reste de $Sh_{p',q'}(\alpha ^{p'}\beta
^{q'})$ quand on ne retient dans chaque shuffle que les $\beta $ de gauche et
les $\alpha $ de droite, et de m\^eme $sh(\beta)$,
$sh(\alpha)$, $sh(\emptyset)$ quand on ne retient que les $\beta
$ de gauche, ou les $\alpha $ de droite, ou aucun $\alpha $ ni $\beta $.\\
Apr\`es gommage de lettres inertes, l'expression $\alpha^i\mu ^k\beta
Sh_{p',q'}(\alpha ^{p'}\beta ^{q'})\alpha\mu ^{m-k}\beta ^j$ devient, selons les
cas~:
$$\begin{matrix}1.1&(k=i=j=0, m=0)\hfill\beta sh(\beta ,\alpha
)\alpha=Y^2sh(\beta ,\alpha ) \\
1.2&(k=i=0,0<j<r,m=0)\hfill\beta sh(\beta )=Ysh(\beta )\\
1.3&(k=i=0, j=r, m>0)\hfill\beta sh(\beta ,\alpha )\alpha \mu ^m\beta
^r=X^mY ^{m+r+2}sh(\beta ,\alpha )\cr
2.1&(k=0,0<i<p,j=0,m=0)\hfill sh(\alpha )\alpha =Ysh(\alpha )\\
2.2&(k=0,0<i<p,0<j<r,m=0)\hfill sh(\emptyset)\cr
2.3&(k=0,0<i<p,j=r,m>0)\hfill sh(\alpha )\alpha \mu ^m\beta ^r=X^mY
^{m+r+1}sh(\alpha )\\
3.1&(k=m,j=0,i=p,m>0)\hfill\alpha ^p\mu ^m\beta sh(\beta
,\alpha )\alpha =X^mY^{m+p+2}sh(\beta ,\alpha )\cr
3.2&(k=m,i=p,0<j<r,m>0)\hfill\alpha ^p\mu ^m\beta sh(\beta)=X^mY
^{m+p+1}sh(\beta )\\
3.3&(0<k<m,i=p,j=r)\hfill\alpha ^p\mu ^k\beta sh(\beta,\alpha
)\alpha \mu ^{m-k}\beta ^r=X^mY^{m+p+r+2}sh(\beta ,\alpha ).\end{matrix}$$

(Remarquons au passage -- bien que cela ne soit pas utile dans la suite
-- qu'un bloc vivant est donc remplac\'e, lors de sa dislocation par $\varphi
B$, par une somme de concat\'enations de $0$, $1$, ou $2$ bloc(s) vivant(s),
intercal\'es \'eventuellements avec des blocs inertes). En utilisant les
coefficients bin\^omiaux, on peut majorer grossi\`erement les $sh$ par~:
$$\displaylines{ sh(\beta ,\alpha)\leq
\sum_{0\leq u\leq p',0\leq v\leq q'}\beta
^v\alpha^uC_{r'-u-v}^{q'-v}= \sum_wY^w
\sum_{0\leq w-v\leq r'-q',0\leq v\leq q'}C_{r'-w}^{q'-v}=\cr
\sum_wY^w
\sum_{w-p'\leq v\leq w,
0\leq v\leq q'}C_{r'-w}^{q'-v}\leq
\sum_{0\leq w\leq r'}Y^w 2^{r'-w},\cr
sh(\beta)\leq\sum_{0\leq v\leq
q'}Y ^vC_{r'-v}^{q'-v},\hfill
sh(\alpha )\leq\sum_{0\leq u\leq p'}Y ^uC_{r'-u}^{p'-u},\hfill
sh(\emptyset)=C_{r'}^{q'}.\cr}$$
Dans chacun des 9 cas, en sommant sur les
valeurs possibles de $i,j,k,p',q'$ (avec $p'+q'=r'=p+q+r-i-j-m$, et
$p-i<p'\leq p+q-m-i$ ou, ce qui revient au m\^eme, $r-j\leq q'<r+q-m-j$), le
nombre $T(m,n)$ de termes de bidegr\'e $(m,n)$ obtenus est donc major\'e par~:
$$\displaylines{(1.1)\hfill{\bf 1}_{(m=0)}\hfill
\sum_{r\leq q'<r+q, 0\leq w\leq p+q+r, n=w+2}\hfill 2^{p+q+r-w}=\cr
\hfill{\bf 1}_{\begin{pmatrix}m=0, & 2\leq n\leq p+q+r+2\end{pmatrix}}\hfill
q2^{p+q+r-n+2},\hfill\cr
(1.2)\hfill{\bf 1}_{(m=0)}\hfill
\sum_{0<j<r, r-j\leq q'<r+q-j, 0\leq v\leq q', n=v+1}\hfill
C_{p+q+r-j-v}^{q'-v}=\cr
\hfill{\bf 1}_{(m=0,0<n)}\hfill
\sum_{0<j<r, r-j\leq q'<r+q-j, n-1\leq q'}\hfill
C_{p+q+r-j-n+1}^{q'-n+1}\leq\cr
\hfill{\bf 1}_{(m=0,0<n)}\hfill
\sum_{0<j<r, j\leq r+q-n}\hfill 2^{p+q+r-j-n+1}\leq\cr
\hfill{\bf 1}_{\begin{pmatrix}m=0,&0<n<r+q, &r>1\end{pmatrix}}\hfill
2^{p+q+r-n+1},\hfill\cr
(1.3)\hfill{\bf 1}_{(m>0)}\hfill
\sum_{0\leq q'<q-m, 0\leq w\leq p+q-m, n-m=w+r+2}
\hfill 2^{p+q-m-w}=\cr
\hfill{\bf 1}_{\begin{pmatrix}0<m<q,&m+r+1<n\leq p+q+r+2\end{pmatrix}}
\hfill(q-m)2^{p+q+r-n+2},\hfill\cr
(2.2)\hfill{\bf 1}_{(m=n=0)}\hfill
\sum_{ 0<i<p,0<j<r,r-j\leq q'<r+q-j }\hfill
C_{p+q+r-i-j}^{q'}\leq\cr
\hfill{\bf 1}_{(m=n=0)}\hfill
\sum_{ 0<i<p,0<j<r }\hfill 2^{p+q+r-i-j}\leq\cr
\hfill{\bf 1}_{\begin{pmatrix}m=n=0,&p>1,&r>1\end{pmatrix}}\hfill 2^{p+q+r},\hfill\cr
(2.3)\hfill{\bf 1}_{(m>0)}\hfill
\sum_{ 0<i<p, p-i<p'\leq p+q-m-i,0\leq u\leq p',n-m=u+r+1 }\hfill
C_{p+q-m-i-u}^{p'-u}=\cr
\hfill{\bf 1}_{(m>0,r<n-m)}\hfill
\sum_{ 0<i<p, p-i<p'\leq p+q-m-i,n-m-r-1\leq p'}\hfill
C_{p+q+r-n-i+1}^{p'+m-n+r+1}\leq\cr
\hfill{\bf 1}_{(0<m<q,r<n-m)}\hfill
\sum_{ 0<i<p, i\leq p+q+r-n+1 }\hfill 2^{p+q+r-i-n+1}\leq\cr
\hfill{\bf 1}_{\begin{pmatrix}0<m<q,&p>1,&m+r<n\leq p+q+r\end{pmatrix}}
\hfill 2^{p+q+r-n+1},\hfill\cr
(3.3)\hfill\sum_{0<k<m,0\leq q'<q-m,0\leq w\leq q-m,n-m=w+p+r+2}
\hfill 2^{q-m-w}=\cr
\hfill{\bf 1}_{\begin{pmatrix}1<m<q,&p+m+r+1<n\leq
p+q+r+2\end{pmatrix}}\hfill(m-1)(q-m)2^{p+q+r-n+2}\hfill\cr}$$ et de m\^eme,
$$\displaylines{(2.1)\hfill{\bf 1}_{\begin{pmatrix}m=0,&0<n\leq
p+q,&p>1\end{pmatrix}}\hfill 2^{p+q+r-n}\hfill\cr
(3.1)\hfill{\bf 1}_{\begin{pmatrix}0<m<q,&p+m+1<n\leq p+q+r+2\end{pmatrix}}
\hfill (q-m)2^{p+q+r-n+2}\hfill\cr
(3.2)\hfill{\bf
1}_{\begin{pmatrix}0<m<q,&r>1,&m+p<n<p+q+r\end{pmatrix}}
\hfill 2^{p+q+r-n+1}.\hfill\cr}$$
En r\'esum\'e,
$$(\varphi B)(X^qY^s)\leq\sum_{0\leq m\leq n, m<q,n\leq s+2}T(m,n)X
^mY^n,$$ avec, en notant $x=q-m$, $y=s-n$~:
$$\begin{matrix}
T(0,n)&\leq &(4q+3)2^y&&\\
T(m,n)&\leq &(1+xq)2^{y+2}&{\rm pour}& m>0\end{matrix}$$
donc dans les deux cas, $T(m,n)\leq q(1+x)2^{y+2} \leq q2^{y+2+{x+1\over 2}}$,
d'o\`u le r\'esultat en posant
$K=2^{5/2},L=2^{1/2},M=2$.\cqfd

\noindent{\sl Preuve du th\'eor\`eme VI.1.\/} Jusqu'\`a pr\'esent, $X^qY
^s$ d\'esignait (le mon\^ome associ\'e \`a) l'un des blocs vivants du mot \`a
transformer par $\varphi B$. Consid\'erons \`a pr\'esent un mot comportant $k$
blocs vivants (associ\'es \`a) $X^{q_i}Y^{s_i}$ ($i=1,\ldots,k$), et
notons ce mot $X^qY^s$ (avec $0\leq q:=\sum q_i\leq s:=\sum s_i$), et
$(\varphi B)(X^qY^s)$ son image. D'apr\`es le lemme,
$$\displaylines{(\varphi B)(X^qY^s)\leq
\sum_{i=1}^k X^{q-q_i}Y^{s-s_i} Kq_i\sum_{0\leq m\leq n, m<q_i,n\leq s_i+2}
L^{q_i}M^{s_i}({X\over L})^m({Y\over M})^n=\cr
KL^qM^s\sum_{i=1}^kq_i\sum_{0\leq m\leq n, m<q_i,n\leq s_i+2}({X\over
L})^{m+q-q_i}({Y\over M})^{n+s-s_i}\leq\cr
KL^qM^s\sum_{i=1}^k q_i\sum_{0\leq u\leq
v,u<q,v\leq s+2}({X\over L})^u({Y\over M})^v=\cr
KL^qM^sq\lfloor{({X\over L})^{q-1}({Y\over M})^{s+2}P({X\over L})P({Y\over
M})}\rfloor, \cr}$$ en posant $P(Z)=1+Z^{-1}+Z^{-2}+\ldots$, et en d\'esignant
par $\lfloor\ \rfloor$ la partie de bidegr\'es $(u,v)$ tels que $0\leq u\leq
v$. On peut donc, en \og~oubliant\fg\ les $\lfloor\ \rfloor$ interm\'ediaires,
majorer les it\'er\'ees de $\varphi B$ par~: $$(\varphi B)^k(X^qY^s)\leq
K^kL^qM^s{q!\over (q-k)!} \lfloor{({X\over L})^{q-k}({Y\over
M})^{s+2k}P^k({X\over L})P^k({Y\over M})}\rfloor.$$ Or pour $k\geq 1$,
$$P^k(Z)=\sum_{i=0}^\infty C_{i+k-1}^{k-1}Z^{-i}\leq \sum_{i=0}^\infty
2^{i+k-1}Z^{-i}=2^{k-1}\sum_{i=0}^\infty({2\over Z})^i.$$ D'o\`u
$$\displaylines{(\varphi B)^k(X^qY^s)\leq
K^kL^qM^s{q!\over (q-k)!}
\sum_{i=0}^{q-k}
({X\over L})^{q-k-i}\sum_{j=0}^{s-q+3k+i}({Y\over
M})^{s+2k-j}2^{i+j+2k-2}=\cr {X^qY^s\over 4}({4KLY^2\over M^2X})^k{q!\over
(q-k)!} \sum_{i=0}^{q-k}({2L\over X})^i\sum_{j=0}^{s-q+3k+i}({2M\over Y
})^j.\cr}$$ En particulier pour un $\omega_n\in\Omega^n(A\otimes B)$
g\'en\'erique, c'est-\`a-dire de la forme\break $a_0\otimes b_0{\rm
d}(a_1\otimes b_1)\ldots{\rm d}(a_n\otimes b_n)$ avec $a_i\in A, b_i\in B$
quelconques, tous de norme 1, on applique la formule avec $1\leq k\leq
q=s=n+1$~: on peut toujours supposer (c'est d'ailleurs le cas pour les
constantes trouv\'ees dans la preuve du lemme) que $2M>1$ et $4LM>1$. Alors,
$(\varphi B)^k(\omega _n)$ est une somme de termes de $\Omega^{n+2k}(A\otimes
B)$, de normes $\leq 1$, dont le nombre est born\'e par~:
$$\displaylines{{1\over 4}({4KL\over M^2})^k{q!\over
(q-k)!} \sum_{i=0}^{q-k}(2L)^i\sum_{j=0}^{3k+i}(2M)^j\cr\leq
{M\over 2(2M-1)}(32KLM)^k{q!\over
(q-k)!} \sum_{i=0}^{q-k}(4LM)^i\cr \leq
{M\over 2(2M-1)(4LM-1)}(4LM)^{q+1}(8K)^k{q!\over(q-k)!}\cr =
{8L^2M^3\over (2M-1)(4LM-1)}(4LM)^n(8K)^k{(n+1)!\over(n+1-k)!}\cr \leq
{8L^2M^3\over (2M-1)(4LM-1)}{(4LM)^n\over\sqrt{n!}}(8K)^k\sqrt{(n+2k)!}.\cr}$$
La norme de $\sum_{k>0}(\varphi B)^k(\omega _n)$
dans $\Omega_r(A\otimes B)$ est donc major\'ee par~:
$$\displaylines{{8L^2M^3\over (2M-1)(4LM-1)}{(4LMr)^n\over\sqrt{n!}}
\sum_{k=1}^{n+1}(8Kr^2)^k\cr\leq   
{64KL^2M^3r^2\over (2M-1)(4LM-1)(1-8Kr^2)}{(4LMr)^n\over\sqrt{n!}},\cr}$$
pourvu que $8Kr^2<1$. Pour tout $r'>0$, $\sum_{k>0}(\varphi B)^k$ d\'efinit
donc bien, pour $r$ assez petit (en fait~: pour $r$ tel que $4LMr\leq r'$ et
$8Kr^2<1$), une application continue de $\Omega_{r'}(A\otimes B)$ dans
$\Omega_r(A\otimes B)$ (et m\^eme~: quand $r\to 0$, la norme de cette
application tend vers $0$).\cqfd

\begin{coro}
$\bbar{sh}_\varepsilon $ est
associatif \`a homotopie pr\`es et $\bbar{AW}_\varepsilon $ est
coassociatif \`a homotopie pr\`es.
\end{coro}

\noindent{\sl Preuve.\/} L'associativit\'e (\`a homotopie pr\`es) de $\bbar{sh}_\varepsilon $
se d\'emontre comme dans la proposition \ref{prop:hassoc} (on v\'erifie facilement que
l'op\'erateur $B_3$ de Getzler-Jones est continu, de $\Omega_\varepsilon
(A)\otimes\Omega_\varepsilon (B)\otimes \Omega_\varepsilon (C)$ dans
$\Omega_\varepsilon (A\otimes B\otimes C)$). La coassociativit\'e \`a
homotopie pr\`es de $AW_\varepsilon $ s'en d\'eduit gr\^ace au th\'eor\`eme
ci-dessus.\cqfd

\end{document}